\documentclass[preprint,12pt]{elsarticle}

\usepackage{graphicx}
\graphicspath{{figures//}}
\usepackage{lineno}
\usepackage{amsmath,amssymb,amsfonts,amsthm}
\usepackage{mathtools}
\usepackage{cases}

\usepackage{caption}
\usepackage{subfig} 
\usepackage{color}

\usepackage{hhline}
\usepackage{multirow}
\usepackage{array}
\usepackage{accents}
\usepackage{booktabs}

\usepackage{hyperref}
\hypersetup{
    colorlinks=true,
    linkcolor=black,
    filecolor=magenta,      
    urlcolor=red,
    pdftitle={Overleaf Example},
    pdfpagemode=FullScreen,
    }

\newcommand{\etp}{e_{t+}}
\newcommand{\etm}{e_{t-}}

\newcommand{\dtp}{\delta_{t+}}
\newcommand{\dtm}{\delta_{t-}}
\newcommand{\dtd}{\delta_{t\cdot}}
\newcommand{\dtt}{\delta_{tt}}

\newcommand{\mtp}{\mu_{t+}}
\newcommand{\mtm}{\mu_{t-}}
\newcommand{\mtd}{\mu_{t\cdot}}

\newcommand{\norm}[1]{\left\|#1\right\|}
\newcommand{\innp}[1]{\left\langle#1\right\rangle}

\usepackage[dvipsnames]{xcolor}
\def\SBcomment[#1]{\textcolor{Red}{#1}}
\def\MDcomment[#1]{\textcolor{Green}{#1}}

\journal{-}

\begin{document}

\begin{frontmatter}


\title{Simulation of the Geometrically Exact Nonlinear String via Energy Quadratisation}



\author[AH]{Michele Ducceschi}
\address[AH]{Department of Industrial Engineering, University of Bologna, IT}
\author[BH]{Stefan Bilbao}
\address[BH]{Acoustics and Audio Group/Reid School of Music, Alison House, University of Edinburgh, UK}

\begin{abstract}
String vibration represents an active field of research in acoustics. Small-amplitude vibration is often assumed, leading to simplified physical models that can be simulated efficiently. However, the inclusion of nonlinear phenomena due to larger string stretchings is necessary to capture important features, and efficient numerical algorithms are currently lacking in this context. Of the available techniques, many lead to schemes which may only be solved iteratively, resulting in high computational cost, and the additional concerns of existence and uniqueness of solutions. Slow and fast waves are present concurrently in the transverse and longitudinal directions of motion, adding further complications concerning  numerical dispersion. This work presents a linearly-implicit scheme for the simulation of the geometrically exact nonlinear string model. The scheme conserves a numerical energy, expressed as the sum of quadratic terms only, and including an auxiliary state variable yielding the nonlinear effects. This scheme allows to treat the transverse and longitudinal waves separately, using a mixed finite difference/modal scheme for the two directions of motion, thus allowing to accurately resolve the wave speeds at reference sample rates. Numerical experiments are presented throughout. 
\end{abstract}


\end{frontmatter}


\section{Introduction}\label{sec:Intro}

The vibration of a musical string, to a first approximation, can be considered linear. From a numerical standpoint, the problem of linear string vibration was explored in many works, see e.g., Ruiz \cite{ruiz_thesis,ruiz_jaes_1971}, Bacon and Bowsher \cite{bacon_bowsher_acustica_1978}, Chaigne and Askenfelt \cite{chaigne_jasa_1994}, Bensa \emph{et. al} \cite{bensa_jasa_2003}, Ducasse \cite{ducasse_jasa_2005}, and others. 

Nonlinear behaviour results for moderate and high amplitudes of motion, and various physical models are available. One of the simplest such models includes tension-modulation effects, as per the Kirchhoff-Carrier model \cite{carrier_QAM_1945}. Numerically, this  was approached using finite difference schemes \cite{bilbao_book} (Chapter 8) as well as Volterra series \cite{helie_JSV_2008}. Longitudinal motion is neglected in this model, which was shown to be inaccurate, since an energy transfer always exists from the transverse to the longitudinal waves \cite{narashima_JSV_1968,rowland_EJP_2011}. 

Finer models, comprising longitudinal motion, can be obtained via geometric arguments. Here, one considers the strain of of a stretched string, and applies Hooke's law. Approximating the energy in a power series results in a model as per Morse and Ingard's \cite{morse_book} (Section 14.3). Numerically, this approximate model was solved via finite difference schemes \cite{bilbao_jasa_2005,bilbao_book,ducceschi_st_ICA_2016}, as well as modal methods, though with further simplifications on the coupling mechanism \cite{bank_JASA_2005}.

In the Morse and Ingard model, the hypothesis of weak linear degenerancy, as explained in \cite{tatsien_book}, is not fulfilled, and this has lead researchers to study the the geometrically exact model. The finite element method, in particular, is preeminent here, see e.g. \cite{chabassier_CMAME_2010,phdthesisChab,marazzato_CMAME_2019}. In \cite{chabassier_CMAME_2010,phdthesisChab}, the conservation of a non-negative numerical energy allows to derive a stability condition involving the mesh size and the time step. Since these methods are fully-implicit,  they require the solution of a large nonlinear algebraic system at each time step, for which  existence and uniqueness of the solution must be proven, and which are only approachable using iterative routines such as Newton-Raphson \cite{chatziioannou_DAFX_2017}. Further complications arise in the choice of the tolerance thresholds and maximum number of iterations of the routine.

For these reasons, in this work, a method is presented, such that the resulting numerical scheme is linearly-implicit in character, thus sidestepping the machinery of iterative methods. Here, existence and uniqueneness of the numerical update are proven by simple inspection of the update matrix, and efficiency is also greatly improved, since the update requires the solution of a single linear system. This particular form is obtained after quadratisation of the nonlinear potential, in a way that is analogous to the Invariant Energy Quadratisation method, proposed orginally for parabolic phase-field models \cite{yang_JCP_2016,yang_CMAME_2017,zhao_IJNM_2017}.
 
Furthermore, the scheme allows for an accurate wideband resolution of both the transverse and longitudinal wave speeds, via a mixed finite difference/modal method employing free parameters. Rather than increasing the formal order of accuracy, the free parameters are used to increase the spectral accuracy across the frequency axis.

The article is structured as follows. In Section \ref{sec:NlinWE} the geometrically nonlinear string is discussed, including stiffness effects. A semi-discrete formulation is given in Section \ref{sec:semidisc}. The fully-discrete numerical scheme is detailed in Section \ref{sec:TimeDisc}, including a proof of conservation of the numerical energy. This section also includes a discussion on numerical dispersion reduction via the use of free parameters. In Section \ref{sec:FullString}, an application comprising the nonlinear string with source and loss terms is illustrated.

\section{Nonlinear Wave Equation}\label{sec:NlinWE}

The case of a nonlinear wave equation describing geometric nonlinearities in stiff strings will be considered here. The geometrically exact nonlinear potential has been used in previous works to simulate the behaviour of musical strings undergoing large amplitude vibration, and comprising a coupling between the longitudinal and transverse components of motion \cite{chabassier_CMAME_2010,chabassier_JASA_2013,marazzato_CMAME_2019}. In the lossless, zero-input case, the system considered here is 
\begin{subequations}\label{eq:VectorString1}
\begin{align}
\left(\rho A\partial_t^2  - T_0 \partial_x^2 + EI \partial_x^4\right) u(t,x) &= \partial_x\left( \frac{\partial \phi}{\partial ({\partial_x u})} \right), \label{eq:VectorString1a} \\
\left(\rho A\partial_t^2  - T_0 \partial_x^2  \right) v(t,x) &= \partial_x\left( \frac{\partial \phi}{\partial ({\partial_x v})} \right), \label{eq:VectorString1b}
\end{align}
\end{subequations}
where 
\begin{equation}\label{eq:PhiContDef}
\phi = \phi(\partial_x u,\partial_x v) = \frac{EA-T_0}{2}\left(\sqrt{(1+\partial_x v)^2+(\partial_x u)^2}-1\right)^2 .
\end{equation}  
First, in \eqref{eq:VectorString1} and in the following, the symbol $\partial^i_j$ denotes a partial derivative along $j$ of order $i$.
The displacements are $u(t,x): \mathbb{R}^+_0 \times [0,L] \rightarrow \mathbb{R}$ in the vertical or flexural direction (perpendicular to the string length), and $v(t,x): \mathbb{R}^+_0 \times [0,L] \rightarrow \mathbb{R}$ in the longitudinal direction (parallel to the string length).

Constants appear as: volume density $\rho$, area of the cross section $A$ (= $\pi r^2$ for a circular cross section of radius $r$), Young's modulus $E$, applied tension $T_0$, string length $L$, area moment of inertia $I$ ($= \pi r^4 / 4$ for a circular cross section). Model  \eqref{eq:VectorString1} is geometrically exact \cite{morse_book,chabassier_CMAME_2010}, and includes a stiffness term as per the Euler-Bernoulli beam model \cite{han_JSV_1999,ducceschi_JASA_2016}. Equations \eqref{eq:VectorString1} must be completed by suitable initial conditions, assumed here differentiable, expressed as 
\begin{equation}\label{eq:InCond}
{ u}(0,x) =  u_0(x), \, { v}(0,x) =  v_0(x), \, {\partial_t u}(0,x) = p_{0}(x), \,{\partial_t v}(0,x) = q_0(x).
\end{equation}

The nonlinear potential $\phi$, as per \eqref{eq:PhiContDef}, is non-negative when $EA \geq T_0$ (a condition that is always satisfied for musical strings). This condition will be assumed valid in the remainder. It is noted that the equations of motion for the geometrically exact nonlinear string are usually expressed in a different form (see e.g. \cite{chabassier_CMAME_2010,marazzato_CMAME_2019}), employing an expression of the potential energy which is non-negative in all cases, though yielding equations identical to \eqref{eq:VectorString1}. The form proposed here, in contrast to the usual form, allows to isolate linear terms proportional to $T_0$, as can be seen in both \eqref{eq:VectorString1a} and \eqref{eq:VectorString1b}. This design choice is justified in a musical context, where strings operate in a regime  driven by the flexural tension term. Here, the linear part of \eqref{eq:VectorString1a} is unaffected by $\phi$. Ultimately, this allows to achieve full audio bandwidth using reference sample rates, as will be shown in Section \ref{subsec:NumDisp}.

\subsection{Energy Identities}\label{subsec:EnIds}
For two square-integrable  functions ${ f}$, ${ g} : [0,L] \rightarrow \mathbb{R}$ one may define an inner product and associated norm as
\begin{equation}\label{eq:innpCnt}
\qquad \innp{{ f},{ g}} = \int_{0}^L { f}  { g} \, \mathrm{d}x, \qquad \norm{{ f}} = \sqrt{\innp{{ f}, { f}}}.
\end{equation}
These  definitions can be used to derive suitable energy identities for the nonlinear wave equation. Here, one takes an inner product of \eqref{eq:VectorString1a} with $\partial_t u$, and of \eqref{eq:VectorString1b} with $\partial_t v$, to get
\begin{subequations}\label{eq:tmpHam}
\begin{align}
\rho A \innp{ \partial_t u, \partial^2_t u} &= T_0 \innp{ \partial_t u, \partial^2_x u} - EI \innp{ \partial_t u, \partial^4_x u}+ \innp{\partial_t u,  \partial_x \left( \frac{\partial \phi}{\partial ({\partial_x u})} \right)}, \\
\rho A \innp{ \partial_t v, \partial^2_t v} &= T_0 \innp{\partial_t v, \partial^2_x v} + \innp{\partial_t v,  \partial_x \left( \frac{\partial \phi}{\partial ({\partial_x v})} \right)}.
\end{align}
\end{subequations}
Summing the two equations, after appropriate integration by parts, and using the fact that
\begin{equation}\label{eq:chnRlPhi}
\frac{d \phi}{dt} =  { \partial_t\partial_{x} u} \left( \frac{\partial \phi}{\partial ({\partial_x u})} \right)  +  {  \partial_t\partial_{x} v}  \left( \frac{\partial \phi}{\partial ({\partial_x v})} \right),
\end{equation}
one gets the following energy balance
\begin{align}
&\frac{dH}{dt} = \left(\innp{\partial_t u, F_u} - \innp{\partial_{t}\partial_x u, M_u} + \innp{\partial_t v, F_v }\right) \big|_0^L,\label{eq:tmpHam2}
\end{align}
where the boundary forces and  moment, due to linear and nonlinear effects, are expressed as
\begin{align*}
F_u = T_0 \partial_x u - EI\partial^3_{x}u +  \frac{\partial \phi}{\partial ({\partial_x u})},\,\,\,
F_v = T_0 \partial_x v +   \frac{\partial \phi}{\partial ({\partial_x v})}, \,\,\, M_u = -EI \partial^2_x u.
\end{align*}
The total energy is the sum of kinetic, linear and nonlinear potential energies, as
\begin{equation}\label{eq:HamWEPhi}
H(t) = E_k(t) + E_{pl}(t) + E_{pnl}(t),
\end{equation}
where
\begin{subequations}\label{eq:HamWEPhi}
\begin{align}
E_k &= \frac{\rho A}{2}\left(\norm{{ \partial_t u}}^2 + \norm{{ \partial_t v}}^2\right), \\
E_{pl} &= \frac{T_0}{2} \left( \norm{\partial_x u}^2 + \norm{\partial_x v}^2  \right) + \frac{EI}{2}  \norm{\partial^2_x u}^2, \\
E_{pnl} & = \norm{\sqrt{\phi}}^2.
\end{align}
\end{subequations}
which is  non-negative. The right hand side of \eqref{eq:tmpHam2} consists of boundary terms. In order for the system to be conservative, it is sufficient to impose conditions such that the boundary terms vanish. Various combinations are possible, and one must be wary of  the nonlinear terms when imposing the vanishing of the forces at the boundary. Here, the string is assumed to be simply-supported flexurally, and fixed longitudinally, i.e.
\begin{equation}\label{eq:BndCnts}
u = \partial_x^2 u = v = 0 \text{ at } x = 0,L
\end{equation}
When these conditions are enforced, energy conservation holds, and thus
\begin{equation}
H(t) = H(0) \triangleq H_0.
\end{equation}
An expression for $H_0$ is derived by direct subsitution of \eqref{eq:InCond} (plus appropriate derivatives) into \eqref{eq:HamWEPhi}. Bounds on solution growth may be expressed easily from such conserved energy, as 
\begin{equation}\label{eq:BoundsWE}
0 \leq \norm{ \partial_t u } \leq \sqrt{2 H_0/ \rho A}, \quad 0 \leq \norm{ \partial_x u} \leq \sqrt{2 H_0/ T_0},
\end{equation}
with  analogous  bounds holding for the longitudinal term. Since the string is fixed at both ends, bounds on $u$, rather than its derivatives, may be obtained here, see e.g. \cite{bilbao_book} (Chapter 6) and \cite{chabassier_CMAME_2010}.

\subsection{Quadratisation}

In view of the numerical application presented below, the nonlinear potential energy is now quadratised by means of the function $\psi : \mathbb{R}^2 \rightarrow \mathbb{R}^+_0$ defined as
\begin{equation}\label{eq:psiDef}
\psi(\partial_x u,\partial_x v) = \sqrt{2 {\phi(\partial_x u,\partial_x v)}}.
\end{equation}
Under such transformation, the total energy is 
\begin{equation}\label{eq:HamWEPsi}
H(t) = E_k(t) + E_{pl}(t) + E_{pnl}(t),
\end{equation}
where  $E_k,E_{pl}$ are the same as \eqref{eq:HamWEPhi}, but where
\begin{equation}
E_{pnl} = \frac{\norm{\psi}^2}{2},
\end{equation}
The total energy here includes only quadratic terms, and the equations of motion read
\begin{subequations}\label{eq:VectorStringPsi}
\begin{align}
\left(\rho A\partial_t^2  - T_0 \partial_x^2 + EI \partial_x^4\right) u(t,x) = \partial_x\left( g_u \psi \right), \label{eq:VectorStringPsi1} \\
\left(\rho A\partial_t^2  - T_0 \partial_x^2  \right) v(t,x) = \partial_x\left( g_v \psi \right), \label{eq:VectorStringPsi2} \\
\frac{d \psi}{dt} = g_u  \partial_t\partial_{x} u + g_v  \partial_t\partial_{x} v.   \label{eq:VectorStringPsi3}
\end{align}
\end{subequations}
These are completed by the identities
\begin{align}
g_u = \frac{\partial \psi}{\partial ({\partial_x u})}, \quad g_v = \frac{\partial \psi}{\partial ({\partial_x v})}. \label{eq:gugv}
\end{align}

The transformation \eqref{eq:psiDef} forms the core of the Invariant Energy Quadratisation method (IEQ), proposed by Yang and associates in the context of nonlinear parabolic phase-field models \cite{yang_JCP_2016,yang_CMAME_2017,zhao_IJNM_2017}, and allowing for a linearly-implicit formulation of the associated numerical schemes.  In parallel, a different kind of quadratisation was proposed by H{\'e}lie, Lopes and Falaize within the context of Port-Hamiltonian systems \cite{lopes_IFAC_2015,phdthesis,falaize:hal-01390501,Lopesphdthesis}. The former method, IEQ, may be applied to any non-negative multivariate potential $\phi$. For the latter method, further requirements on convexity must be met by the function $\phi$  \cite{Lopesphdthesis}, precluding the possibility of treating various cases such as e.g. non-invertible potentials. Furthermore, for the associated numerical schemes, second-order accuracy is achieved via a two-step procedure in the latter method, involving the solution of two linear systems per time step, as opposed to one single step for IEQ. A check of the order of accuracy, in the case of a scalar nonlinear ODE, is offered in \ref{app:DuffErrs}.

IEQ for the cases of nonlinear hyperbolic wave equations were proposed in \cite{jiang_JSC_2019} for the sine-Gordon equation, and in \cite{ducceschi_ICA_2019} for the geometrically exact string, and the results therein are extended here. Other examples of the application of IEQ, of interest in musical acoustics, are given in \cite{ducceschi_jasa_2021}.

\section{Semi-discrete formulation}\label{sec:semidisc}

The quadratised equations \eqref{eq:VectorStringPsi} are now discretised in space, using a mixed formulation including a finite difference discretisation of the transverse waves, and a spectral discretisation of the longitudinal waves. The reason for this design choice resides in the difficulty of resolving a system that comprises two different wave speeds (here, transverse and longitudinal), and whose expressions are obtained after linearisation of  \eqref{eq:VectorString1}. The two speeds are, respectively, $c_u = \sqrt{T_0/\rho A}$ (valid at low frequencies), and $c_v = \sqrt{E/\rho}$, where typically $c_v \gg c_u$ for musical strings. Since the grid spacing in a finite difference scheme is directly related to the wave speed via a stability condition, one should consider that adapting the grid spacing to match the velocity in either direction of motion will result in severe frequency warping effects in the other direction, at reference sample rates. One possible solution, proposed e.g. in \cite{bilbao_book} (Chapter 8) and \cite{torin_DAFX_2013}, is to make use of two separate grids, one for the longitudinal and one for the transverse waves, and to perform interpolation in order to couple them. Another choice, explored here, is to make use of a mixed approach: the transverse waves may still be resolved via finite difference schemes; the longitudinal waves may be instead approximated by a suitable modal or spectral discretisation \cite{gazdag_JCP_1973,gazdag_geophysics_1981}.

\subsection{Spatial Finite Difference Operators}

For the transverse direction, the length $L$ is divided into $N$ subintervals, yielding $N+1$ grid points including the end points. The subintervals are of length $h = L/N$, the grid spacing. Boundary conditions are of fixed type, as per \eqref{eq:BndCnts}, therefore the two end points need not be computed and stored. The physical displacement $u(t,x=mh)$ is approximated by $({\bf u}(t))_m$, $1\leq m \leq N-1$.  The spatial difference operator, for fixed conditions, can be expressed as a $N \times N-1 $ matrix of the following form
\begin{equation}
{\bf D}^-{\bf u}  = \frac{1}{h}\left([{\bf u}^\intercal,0]^{\intercal} - [0,{\bf u}^\intercal]^{\intercal}\right).
\end{equation}
The forward and backward spatial difference matrices are related by a transpose operation, i.e. 
\begin{equation}\label{eq:DxmTr}
{\bf D}^{+} = -\left({\bf D}^{-}\right)^\intercal.
\end{equation}
Furthermore, one has
\begin{equation}\label{eq:DxxTr}
{\bf D}^{2} = {\bf D}^{+}{\bf D}^{-}, \qquad {\bf D}^{4} =  {\bf D}^{2} {\bf D}^{2}.
\end{equation}
i.e. the non-commutative composition of the backward and forward differentiation yields the second spatial difference, and the composition of the second spatial difference with itself gives the fourth spatial difference.  Both these matrices are of size $N-1 \times N-1$. Using Taylor-expansion arguments, one has
\begin{align}
\left({\bf D}^{(+,-)}{\bf u}\right)_m &= \frac{du(x)}{dx}\big|_{x=mh} + O(h), \\
\left({\bf D}^{2}{\bf u}\right)_m   &= \frac{d^2u(x)}{dx^2}\big|_{x=mh} + O(h^2), \\
\left({\bf D}^{4}{\bf u}\right)_m   &= \frac{d^4u(x)}{dx^4}\big|_{x=mh} + O(h^2).
\end{align}

\subsection{Spectral Discretisation of the Longitudinal Waves}

For the longitudinal component, one may approximate the continuous function $v(t,x)$ by a grid function ${\bf v}(t): {\mathbb R}^{+}_0 \rightarrow \mathbb{R}^{N-1}$. Then, a spectral operator is applied. This is 
\begin{equation}\label{eq:vExpans}
{\bf v}(t) = {\bf Z}{\bf s}(t). 
\end{equation}
Explicitly, this is given as
\begin{equation}\label{eq:vZ}
v_m(t) = \sqrt{\frac{2h}{L}} \sum_{\nu=1}^{N_s}  \sin \left( \frac{m \nu  h \pi}{L}\right) s_\nu(t) \,\,, \quad 1 \leq m \leq N-1,
\end{equation}
Thus, $Z_{m,\nu}=\sqrt{2h/L}\sin(m\nu h\pi /L)$, and the matrix is of size $\left(N-1\right)\times N_s$. Here, the upper bound in the sum, $N_s$, represents the number of eigenfunctions in the longitudinal direction, to be specified later. Note that such basis functions are consistent with the fixed conditions at the string's ends, as per \eqref{eq:BndCnts}. The square root factor multiplying the basis functions is here only a useful normalisation. 
Here, ${\bf Z}$  satisfies the following identities
\begin{equation}\label{eq:specpropr2}
{\bf Z}^\intercal{\bf Z}  = {\bf I}, \quad {\bf Z}^\intercal\,{\bf D}^{2}\,{\bf Z} = -{\bf \Lambda}, 
\end{equation}
where ${\bf \Lambda}$ is an $N_s\times N_s$ diagonal matrix with $[{\bf \Lambda}]_{\nu,\nu} = \frac{\nu^2 \pi^2}{L^2} + O(h^2)$. The first property is  a statement of the orthogonality of sine functions; the second property is a consequence of the spectral decomposition of the ${\bf D}^2$ operator.
These properties can be used to represent the scheme for the longitudinal wave equation in a modal form, as shown below.

\subsection{Auxiliary State Variable}

The continuous function $\psi(\partial_x u,\partial_x v)$ in \eqref{eq:VectorStringPsi} will be approximated here as an extra, independent state variable, defined on an interleaved spatial grid, i.e. ${\boldsymbol \psi}(t): \mathbb{R}^+_0 \rightarrow {\mathbb{R}^N}$. 
Explicit realisations for \eqref{eq:gugv} must also be supplied. These are given via the vectors ${\bf g}_u$, ${\bf g}_v$:
\begin{equation}\label{eq:gdefsWE}
{\bf g}_u = \frac{\sqrt{EA-T_0}\,\,{\bf D}^{-}{\bf u}}{\sqrt{(1+{\bf D}^{-}{\bf v})^2 + ({\bf D}^{-}{\bf u})^2} }, \quad {\bf g}_v = \frac{\sqrt{EA-T_0}\,\,(1+{\bf D}^{-}{\bf v})}{\sqrt{(1+{\bf D}^{-}{\bf v})^2 + ({\bf D}^{-}{\bf u})^2} },
\end{equation} 
where the division of vectors, and exponentiation, are intended elementwise. Both these vectors are of length $N$. From these, the square $N \times N$ matrices ${\bf G}_u$ and ${\bf G}_v$ are given as
\begin{equation}\label{eq:Gs}
{\bf G}_u = \text{diag}({\bf g}_u), \quad {\bf G}_v = \text{diag}({\bf g}_v).
\end{equation}

\subsection{Semi-discrete Equations of Motion}

Using the proposed notation, a semi-discrete realisation of \eqref{eq:VectorStringPsi} is given as
\begin{subequations}\label{eq:SemiDisc}
\begin{align}
\left(\rho A \frac{d^2}{dt^2} - T_0{\bf D}^2 + EI {\bf D}^4\right){\bf u}(t) = {\bf D}^+ {\bf G}_u {\boldsymbol \psi}(t), \label{eq:SemiDisc1} \\
\left(\rho A \frac{d^2}{dt^2} - T_0{\bf D}^2 \right){\bf Z}{\bf s}(t) = {\bf D}^+ {\bf G}_v {\boldsymbol \psi}(t), \label{eq:SemiDisc2} \\
\frac{d{\boldsymbol \psi}(t)}{dt} = {\bf G}_u \left(\frac{d}{dt}{\bf D}^- {\bf u}(t) \right) + {\bf G}_v \left(\frac{d}{dt}{\bf D}^- {\bf Z}{\bf s}(t) \right). \label{eq:SemiDisc3}
\end{align}
\end{subequations}
Explicit expressions for ${\bf G}_u$, ${\bf G}_v$ are given in \eqref{eq:Gs}. The longitudinal displacement is expressed as a superposition of modes, as per \eqref{eq:vExpans}, and diagonalisation of \eqref{eq:SemiDisc2}  can be obtained by multiplying on the left by ${\bf Z}^\intercal$, to give
\begin{equation}
\left(\rho A \frac{d^2}{dt^2} + T_0 {\bf \Lambda} \right){\bf s}(t) = {\bf Z}^\intercal\,{\bf D}^+ {\bf G}_v {\boldsymbol \psi}(t).
\end{equation}
\subsection{Semi-discrete Energy Identities}

A discrete version of \eqref{eq:innpCnt} (inner product and associated norm) for two vectors ${\bf f}, {\bf g}: \mathbb{R}^+_0 \rightarrow \mathbb{R}^{N-1}$  can be given as
\begin{equation}\label{eq:InnpDisc}
\innp{{\bf f}, {\bf g}} = h \,{\bf f}^\intercal \, {\bf g}, \qquad \norm{{\bf f}} = \sqrt{\innp{{\bf f}, {\bf f}}}.
\end{equation}
Thus, taking an inner product of \eqref{eq:SemiDisc1} with $\frac{d{\bf u}}{dt}$, of \eqref{eq:SemiDisc2} with $\frac{d{\bf v}}{dt}$, and summing, yields a semi-discrete energy balance of the kind 
\begin{equation}
\frac{d\mathfrak{h}(t)}{dt} = 0 \,\, \text{ where } \,\, \mathfrak{h}(t) = \mathfrak{E}_k(t) + \mathfrak{E}_{pl}(t) + \mathfrak{E}_{pnl}(t).
\end{equation}
The semi-discrete kinetic, linear and nonlinear potential energies are
\begin{subequations}\label{eq:SemiDiscEnergy}
\begin{align}
\mathfrak{E}_k & = \frac{\rho A}{2}\left(\norm{\frac{d{\bf u}}{dt}}^2 + \norm{\frac{d{\bf v}}{dt}}^2 \right), \\
\mathfrak{E}_{pl} & = \frac{T_0}{2} \left(\norm{{\bf D}^-{\bf u}}^2 + \norm{{\bf D}^-{\bf v}}^2 \right) + \frac{EI}{2}\norm{{\bf D}^2{\bf u}}^2, \\
\mathfrak{E}_{pnl} &= \frac{\norm{{\boldsymbol \psi}}^2}{2}.
\end{align}
\end{subequations}
which is  a non-negative, semi-discrete counterpart of \eqref{eq:HamWEPsi}. In order to obtain this expression for the energy, the transposition and symmetry properties \eqref{eq:DxmTr}, \eqref{eq:DxxTr} of the difference matrices were used, along with the rate of change of the auxiliary state variable $\boldsymbol \psi$, as per \eqref{eq:SemiDisc3}.

\section{Time discretisation}\label{sec:TimeDisc}

A numerical scheme arises after an appropriate time discretisation of \eqref{eq:SemiDisc}. Thus, the grid functions ${\bf u}(t),{\bf v}(t)$ are approximated at the time $nk$ by the time series ${\bf u}^n, {\bf v}^n$, where $n \in \mathbb N_0$, and where $k$ is the time step (the multiplicative inverse of the sample rate). Similarly, the auxiliary state variable ${\boldsymbol \psi}(t)$ is approximated by an interleaved time series ${\boldsymbol \psi}^{n-1/2}$.

\subsection{Time Difference Operators}
The basic operators in discrete time are the identity and shift operators, defined as
\begin{equation}
   {1}{\bf u}^n = {\bf u}^n, \quad \etp {\bf u}^n = {\bf u}^{n+1}, \quad \etm {\bf u}^n = {\bf u}^{n-1}.
\end{equation}
From these, one may define the time difference operators, all approximating the first time derivative, as
\begin{subequations}
\begin{align}
    \dtp {\bf u}^n &= \frac{(\etp - 1){\bf u}^n}{k} = \frac{d{\bf u}(t)}{dt}\big|_{t=kn} + O(k),\\ 
    \dtm {\bf u}^n &= \frac{(1 - \etm){\bf u}^n}{k}= \frac{d{\bf u}(t)}{dt}\big|_{t=kn} + O(k), \\ 
    \dtd {\bf u}^n &= \frac{(\etp - \etm){\bf u}^n}{2k} = \frac{d{\bf u}(t)}{dt}\big|_{t=kn} + O(k^2).
\end{align}
\end{subequations}
An approximation to the second time derivative is constructed from the above as
\begin{equation}
    \dtt u^n = (\dtp\dtm){\bf u}^n = \frac{d^2{\bf u}(t)}{dt^2}\big|_{t=kn} + O(k^2).
\end{equation}
Averaging operators are also used throughout the text, and are
\begin{subequations}
\begin{align}
    \mtp {\bf u}^n &= \frac{(\etp + 1){\bf u}^n}{2} = {\bf u}(kn) + O(k),\\ 
    \mtm {\bf u}^n &= \frac{(1 + \etm){\bf u}^n}{2}= {\bf u}(kn) + O(k), \\ 
    \mtd {\bf u}^n &= \frac{(\etp + \etm){\bf u}^n}{2} = {\bf u}(kn) + O(k^2).
\end{align}
\end{subequations}
 For the interleaved function ${\boldsymbol \psi}^{n-1/2}$, the same definitions apply formally, but the order of the approximation changes. Thus
\begin{equation}
   1{\boldsymbol \psi}^{n-1/2} = {\boldsymbol \psi}^{n-1/2}, \quad \etp {\boldsymbol \psi}^{n-1/2} = {\boldsymbol \psi}^{n+1/2}.
\end{equation}
The time difference is
\begin{equation}
    \dtp {\boldsymbol \psi}^{n-1/2} = \frac{(\etp - 1){\boldsymbol \psi}^{n-1/2}}{k} = \frac{d{\boldsymbol \psi}(t)}{dt}\big|_{t=kn} + O(k^2),
\end{equation}
and the averaging operator gives
\begin{equation}
\mtp {\boldsymbol \psi}^{n-1/2}  = \frac{(\etp + 1){\boldsymbol \psi}^{n-1/2}}{2} = {\boldsymbol \psi}(kn) + O(k^2).
\end{equation}
and note that the following identity, used throughout the text, holds
\begin{equation}\label{eq:psiId}
\mtp {\boldsymbol \psi}^{n-1/2} = \frac{k}{2}\dtp {\boldsymbol \psi}^{n-1/2} + {\boldsymbol \psi}^{n-1/2}
\end{equation}
Three useful identities are given here. Considering the inner product and norm given in \eqref{eq:InnpDisc}, one has
\begin{subequations}\label{eq:IdFD}
\begin{align}
 \innp{\dtd {\bf u}, \dtt {\bf u}}  &= \dtp \frac{\norm{\dtm {\bf u}}^2}{2}, \\ 
 \innp{{\bf u}, \dtd {\bf u}} &= \dtp \frac{\innp{{\bf u},\etm {\bf u}}}{2}, \\ 
 \innp{\dtp{\boldsymbol \psi},\mtp{\boldsymbol \psi}}  &= \dtp\frac{\norm{\boldsymbol \psi}^2}{2}.
 \end{align}
\end{subequations}

\subsection{Fully-discrete Equations of Motion}

Given the above definitions, a fully-discrete system of equations can now be given. These are
\begin{subequations}\label{eq:FullDisc}
\begin{align}
\left(\rho A \dtt - T_0{\bf D}^2 + EI {\bf D}^4\right){\bf u}^n = {\bf D}^+ {\bf G}_u \,\,\mtp{\boldsymbol \psi}^{n-1/2}, \label{eq:FullDisc1} \\
\left(\rho A \dtt + T_0{\bf \Lambda} \right){\bf s}^n = {\bf Z}^\intercal{\bf D}^+ {\bf G}_v \,\,\mtp{\boldsymbol \psi}^{n-1/2}, \label{eq:FullDisc2} \\
\dtp {\boldsymbol \psi}^{n-1/2} = {\bf G}_u \left(\dtd {\bf D}^- {\bf u}^n \right) + {\bf G}_v \left(\dtd {\bf D}^- {\bf Z}{\bf s}^n \right). \label{eq:FullDisc3}
\end{align}
\end{subequations}
which  discretises \eqref{eq:SemiDisc}. It is remarked that this scheme is a three-step scheme, in that not only does one need to solve for the displacement ${\bf u}$, and modal coordinates ${\bf s}$, but also for ${\boldsymbol \psi}$, which is treated here as an auxiliary state variable. Thus, the numerical system is here completed by a set of initial conditions on ${\bf u}, {\bf v}$ discretising \eqref{eq:InCond}, (from which one may get ${\bf s} = {\bf Z}^\intercal {\bf v}$ ), as well as a suitable initial condition for ${\boldsymbol \psi}$, which can be given as
\begin{equation}
{\boldsymbol \psi}^{1/2} = \sqrt{EA-T_0}\left( \sqrt{\left(1 + {\bf D}^-{\bf Z}\mtp{\bf s}^0\right)^2 + \left({\bf D}^-\mtp{\bf u}^0\right)^2} - 1\right).
\end{equation}
Note that the averaging operators are here used so to yield a second-order accurate discretisation of the continuous function  $\psi$, at the initial interleaved time step $t_{n=1} = k/2$.

The time discretisation was chosen here such that \eqref{eq:FullDisc} conserves a numerical energy. This is easily shown after taking an inner product (as per \eqref{eq:InnpDisc}) of \eqref{eq:FullDisc1} with $\dtd {\bf u}$, of \eqref{eq:FullDisc2} with $\dtd {\bf v}$, summing, and making use of \eqref{eq:FullDisc3} plus the identities given in \eqref{eq:IdFD}. The result is
\begin{equation}\label{eq:FullDiscEner}
\dtp {\mathfrak h}^{n-1/2} = 0, \,\, \text{ where } \,\, \mathfrak{h} = \mathfrak{E}_k + \mathfrak{E}_{pl} + \mathfrak{E}_{pnl}.
\end{equation}
The kinetic, linear and nonlinear potential energies are given as
\begin{subequations}\label{eq:DiscEnergy}
\begin{align}
\mathfrak{E}^{n-1/2}_k &= \frac{\rho A}{2}\left(\norm{\dtm {\bf u}^n}^2 + \norm{\dtm {\bf v}^n}^2\right), \label{eq:DiscEnergy1} \\
\mathfrak{E}^{n-1/2}_{pl} &= \frac{T_0}{2}\left(\innp{{\bf D}^-{\bf u}^n,{\bf D}^- \etm {\bf u}^n} + \innp{{\bf D}^-{\bf v}^n,{\bf D}^- \etm {\bf v}^n}\right) + \frac{EI}{2} \innp{{\bf D}^2{\bf u}^n,{\bf D}^2 \etm {\bf u}^n},  \label{eq:DiscEnergy2} \\
\mathfrak{E}^{n-1/2}_{pnl} &= \frac{\norm{\boldsymbol \psi^{n-1/2}}^2}{2}. \label{eq:DiscEnergy3}
\end{align}
\end{subequations}
This expression approximates the total energy of the system. As opposed to the continuous and the semi-discrete cases \eqref{eq:HamWEPsi} and \eqref{eq:SemiDiscEnergy}, the fully-discrete expression for the energy may be negative in some cases. Non-negativity is achieved here only when the grid spacing $h$ and the time step $k$ satisfy a suitable condition. In this regard, it is noted that such condition arises solely as a consequence of the explicit discretisation of the linear part, since the nonlinear energy, from \eqref{eq:DiscEnergy3}, is non-negative by definition. Thus, a check on the eigenvalues of the linear part allows to conclude that a necessary and sufficient condition for the non-negativity of the total energy is given by (see e.g. \cite{bilbao_book,ducceschi_WaveMotion_2019})
\begin{equation}\label{eq:StabCndStiff}
h \geq  \sqrt{\frac{  T_0k^2  + \sqrt{ ( T_0k^2 )^2 + 16 \rho A E I k^2 } }{2\rho A }}.
\end{equation}
It is remarked that this stability condition is not CFL-like, since the time step is not directly proportional to the grid spacing. This is a consequence of the explicit discretisation of the fourth spatial derivative, using the matrix ${\bf D}^4$. Whilst still resolving waves at full bandwidth, this discretisation has the benefit of reducing the number of grid points, thus improving efficiency, compared to an implicit discretisation with a pure CFL condition. This claim is assessed in  \ref{app:ExpVsImp}.

When  condition \eqref{eq:StabCndStiff} is enforced, the energy is non-negative, and hence the grid functions remain bounded over time, with bounds holding here as discrete versions of \eqref{eq:BoundsWE}:
\begin{equation}
0 \leq \norm{\dtm {\bf u}} \leq \sqrt{2{\mathfrak h}^{1/2}/\rho A}, \quad 0 \leq \norm{\mtm {\bf D}^-{\bf u}} \leq \sqrt{2{\mathfrak h}^{1/2}/T_0},
\end{equation}
with similar bounds holding for the longitudinal grid function $\bf v$. It is remarked that the auxiliary state variable ${\boldsymbol \psi}$ remain itself bounded, as 
\begin{equation}
0 \leq \norm{{\boldsymbol \psi}} \leq \sqrt{2{\mathfrak h}^{1/2}}.
\end{equation}

\subsection{Number of Numerical Longitudinal Eigenfunctions }
A bound on the  number of longitudinal eigenfunctions, determining the size of $\bf Z$ in \eqref{eq:vExpans}, can be derived by imposing non-negativity of the longitudinal linear discrete energy, i.e.  $N_s \leq (2L/\pi k)\sqrt{\rho A / T_0}$ \cite{bilbao_book} (Chapter 6). It will be convenient, however, to use a much lower number of modes, such that
\begin{equation}\label{eq:Ns}
N_s \leq \frac{2L}{\pi k}\sqrt{\frac{\rho}{E}}.
\end{equation}
This is the bound on the number of modes associated with the CFL condition on the longitudinal wave speed $c_v = \sqrt{E/\rho}$. This choice has two benefits: the first is efficiency, since the system is much smaller in size (typically $c_v \gg c_v$); second, frequency warping effects can be minimised, as will be shown in Section \ref{subsec:NumDisp}.

\subsection{Solvability, Local Truncation Error and Convergence}

Solving system \eqref{eq:FullDisc} is accomplished by first using identity \eqref{eq:psiId} in both \eqref{eq:FullDisc1} and \eqref{eq:FullDisc2}, and then using \eqref{eq:FullDisc3} to express the $\dtp {\boldsymbol \psi}$ in terms of $ {\bf u}$, $ {\bf s}$. The resulting system is written as
\begin{equation}\label{eq:Update}
\begin{bmatrix}\frac{\rho A}{k^2}{\bf I} - \frac{1}{4}{\bf D}^+ {\bf G}_u^2 {\bf D}^- & - \frac{1}{4}{\bf D}^+ {\bf G}_u {\bf G}_v {\bf D}^-{\bf Z} \\ - \frac{1}{4}{\bf Z}^\intercal{\bf D}^+ {\bf G}_u {\bf G}_v {\bf D}^-  & \frac{\rho A}{k^2}{\bf I} - \frac{1}{4}{\bf Z}^\intercal{\bf D}^+ {\bf G}_v^2 {\bf D}^-{\bf Z}\end{bmatrix}\begin{bmatrix}{\bf u}^{n+1} \\ {\bf s}^{n+1}\end{bmatrix} = {\bf b}^n,
\end{equation}
where ${\bf b}^n$ is a vector of known coefficients from previous time steps. The update matrix is a non-singular, symmetric matrix, and therefore  ${\bf u}^{n+1}$, ${\bf s}^{n+1}$ are uniquely determined. These values can then be used to compute ${\bf v}^{n+1}$ via \eqref{eq:vExpans}, and ${\boldsymbol \psi}^{n+1/2}$ via \eqref{eq:FullDisc3}. The top left square block in the update matrix is sparse, and it contains the most elements in the matrix. The bottom right square block is the smallest, and is full. The  off-diagonal reactangular blocks are also full. It is convenient to approach the solution of the linear system via the Schur complement, though implementation details are not discussed further here.

The local truncation error $\tau^n_m$ of scheme \eqref{eq:FullDisc}  at time $t_n = nk$ and at the grid point $x_m = mh$ is obtained when the scheme is applied to the true solutions $u(t,x), v(t,x)$ of \eqref{eq:VectorStringPsi}, and Taylor-expanding, see e.g. \cite{leveque_book}. Using the Taylor expansions given above for the spatial and temporal operators, it can be  shown  that
\begin{equation}
\tau^n_m = O(h^2) + O(k),
\end{equation}
and thus the local truncation error is second-order  in $h$ and first-order in $k$. Note that first-order accuracy in $k$ is a consequence of the explicit discretisation of the fourth-order differential operator. Since $k$ and $h$ are chosen along the path dictated by the stability condition  \eqref{eq:StabCndStiff}, one has $k = O(h^2)$, and the scheme is also second-order in the pathlength $s$, as shown in \ref{app:thetaErrs} for the linear part.

One may go further and prove convergence as a consequence of stability and consistency, showing that the scheme converges to the true solution as $h$, $k$ tend to zero, and that the global error preserves the same order accuracy when $h$, $k$ are decreased along the path given by the stability condition \eqref{eq:StabCndStiff}. This is a delicate point: one should be wary that theoretical error estimate blow-ups have been predicted, even for linear wave equations, see e.g. \cite{joly_book,chabassier_CRM_2017}. In \cite{joly_book}, it is shown that a post-processing of the error time series leads indeed to an estimate that does not blow up, in accordance to empirical observations. Here, a formal proof for the proposed scheme, while necessary, is out-of-scope. A proof of order-accuracy in the case of a hyperbolic wave equation as per the sine-Gordon model, using IEQ, is given in \cite{jiang_JSC_2019}. Other formal proofs of applications of IEQ are included in e.g. \cite{yang_JCP_2016}. A check on the accuracy of the discretisation of the linear part, for the transverse direction, is given in \ref{app:thetaErrs}.

\subsection{Experiments}\label{sec:experiments}
\begin{figure}[hbt]
\centering
\includegraphics[width=0.95\linewidth,clip, trim = {3cm 1cm 2cm 1cm}]{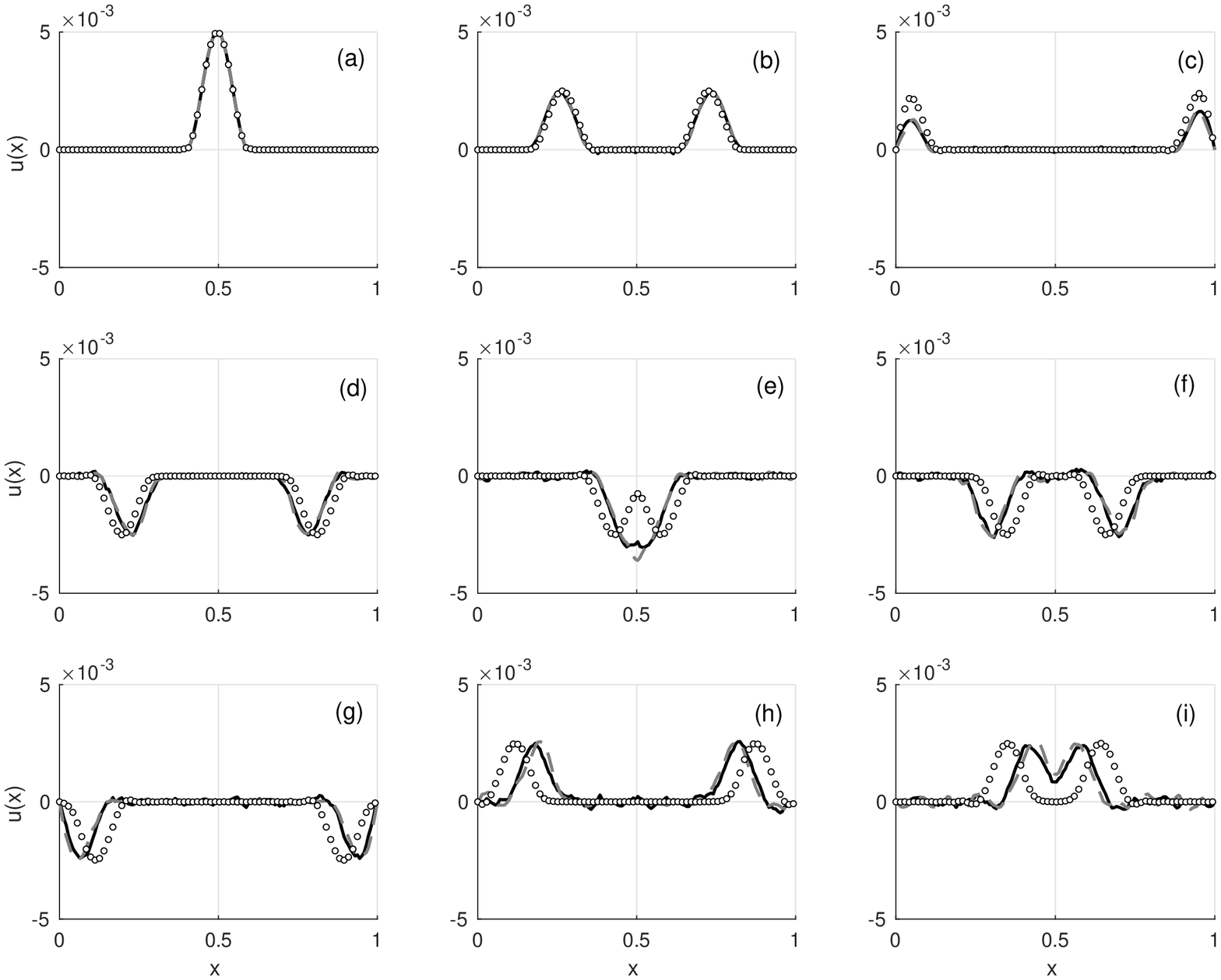}
\caption{Snapshots of the time evolution of solutions to the nonlinear wave equation, at successive times. For all panels, the circles represent the solution to the linear wave equation, the thick black line is the iterative finite difference scheme presented in \cite{ducceschi_ICA_2019}, and the dashed grey line is the scheme  \eqref{eq:FullDisc}. Note that stiffness was not considered in this example, so $I = 0$. The string has parameters $\rho = 8000$ kg/m$^3$, $E = 2\cdot 10^{11}$ Pa, $r = 0.2$ mm, $L = 1$ m, $T_0 = 50$ N. The schemes are initialised with a raised cosine distribution with compact support of the form $u(0,x) = \frac{U_0}{2} \left(1 + \cos \left( \frac{2 \pi (x-L/2)}{2\sigma L}  \right)\right)$ for $L/2-\sigma L \leq x \leq L/2 + \sigma L$, and $u(0,x) = 0$ elsewhere. Here, $U_0 = 5$ mm and $\sigma = 0.1$. The sample rate used is $f_s = 48$ kHz, and the grid spacing is $h = 1.5\sqrt{T_0/\rho A}\,k$}\label{fig:snapsFD}
\end{figure}

\begin{figure}[hbt]
\centering
\includegraphics[width=0.95\linewidth,clip, trim = {1cm 1cm 2cm 1cm}]{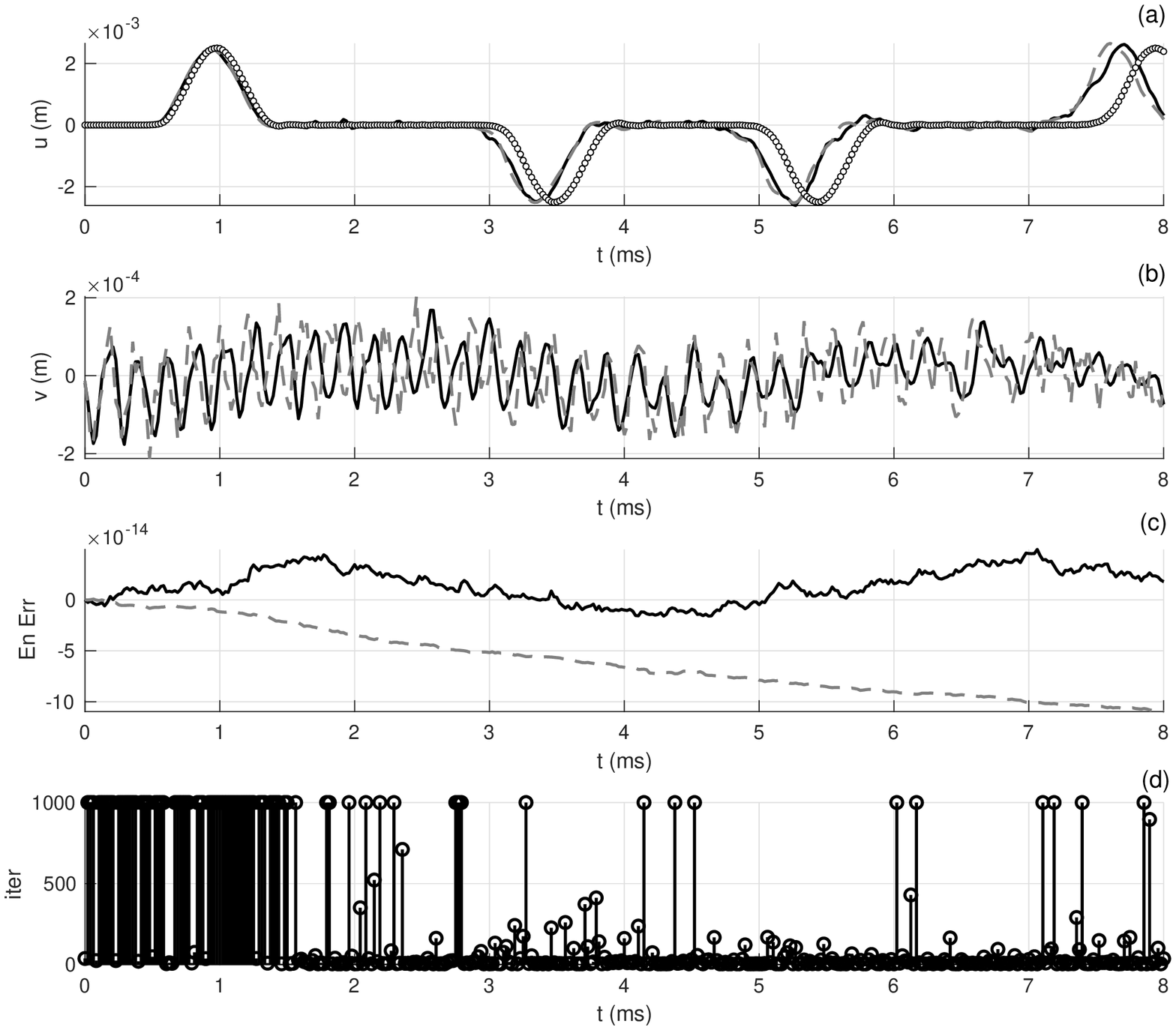}
\caption{Waveforms and energy error. The parameters used are the same as Figure \ref{fig:snapsFD}. The circles represent the solution to the linear wave equation, the thick black line is the finite difference iterative scheme from \cite{ducceschi_ICA_2019}, and the dashed grey line is the non-iterative scheme \eqref{eq:FullDisc}. (a): transverse displacement, output is taken at $x_o=0.72L$; (b): longitudinal displacement, output is taken at $x_o=0.72L$; (c): energy error for both schemes, defined as  $\epsilon^{n-1/2}=1-\mathfrak{h}^{n-1/2}/\mathfrak{h}^{1/2}$; (d): number of iterations for the transverse Newton-Raphson, with a tolerance threshold $\tau=10^{-15}$, and escape condition iter $\leq 1000$. }\label{fig:WaveFormsFD}
\end{figure}

The performance of scheme \eqref{eq:FullDisc} is now compared in a numerical experiment. In Fig. \ref{fig:snapsFD}, one can see the snapshots of the time evolution of the simulation, at times as indicated. The string considered here has physical and geometrical parameters resembling those of a typical musical string. The string is initialised in the transverse direction with a raised cosine distribution centered around the midpoint of the problem domain. In this example, the stiffness term $EI$ was set to zero, so to compare the results against both the linear one-dimensional wave equation, and to the iterative finite-difference scheme presented in \cite{ducceschi_ICA_2019}. The key features of nonlinear wave propagation relative to linear are increased wave speed, as well as a progressive flattening of the peaks with time. Notice that both schemes \eqref{eq:FullDisc} and the iterative scheme from \cite{ducceschi_ICA_2019} yield consistent solutions.
The same observations can be drawn from results shown in Fig. \ref{fig:WaveFormsFD}. The top panel is a representation of the same solutions for the transverse displacement, but plotted against time. The second panel instead shows the longitudinal displacement: the two models yield again consistent results. Notice that, in this case, the overall magnitude of the longitudinal motion is about one order of magnitude less than that of the transverse. Energy remains conserved for both schemes, with an error of the order of machine accuracy. However, the iterative scheme can only achieve this through a large number of iterations, highlighting the advantage of the linearly-implicit formulation \eqref{eq:Update}.

Note that, in the simulation, the grid spacing is chosen to be larger than bound \eqref{eq:StabCndStiff} satisfied with equality. This is to avoid possible erratic behaviour near Nyquist. In practice, some spurious noise is generated by the scheme, and propagated through the spectrum by the nonlinear coupling. This is not a blow-up, rather, it is a polluted solution. This effect has been observed for other nonlinear systems, see for example the Kirchhoff-Carrier string in \cite{bilbao_book} (Chapter 8).

\subsection{Wideband Numerical Dispersion Reduction}\label{subsec:NumDisp}

Numerical dispersion introduces artefacts in the computed solution. In the linear case, modified-equation techniques can be employed to increase the formal order of accuracy of the numerical schemes, see e.g. \cite{cohen_1996,joly_2010,dablian_2010,warming_1974,bilbao_2018}. Dispersion is  a frequency-dependent effect, and higher frequencies are usually more affected by it. Order-accuracy holds in the low-frequency limit, and it was observed that lower-order accurate schemes may yield better wideband behaviour, see e.g. \cite{germain_dafx_2015,ducceschi_DAFX_2021}. For this reason, it may be preferable to reduce dispersion over the entire range of frequencies \cite{lele_1992,ducceschi_WaveMotion_2019}, rather than increasing the formal order of accuracy. The linear parts of \eqref{eq:VectorString1a} and \eqref{eq:VectorString1b} will now be parameterised.

\begin{figure}[hbt]
\centering
\includegraphics[width=0.92\linewidth,clip, trim = {0cm 0.0cm 1cm 0.0cm}]{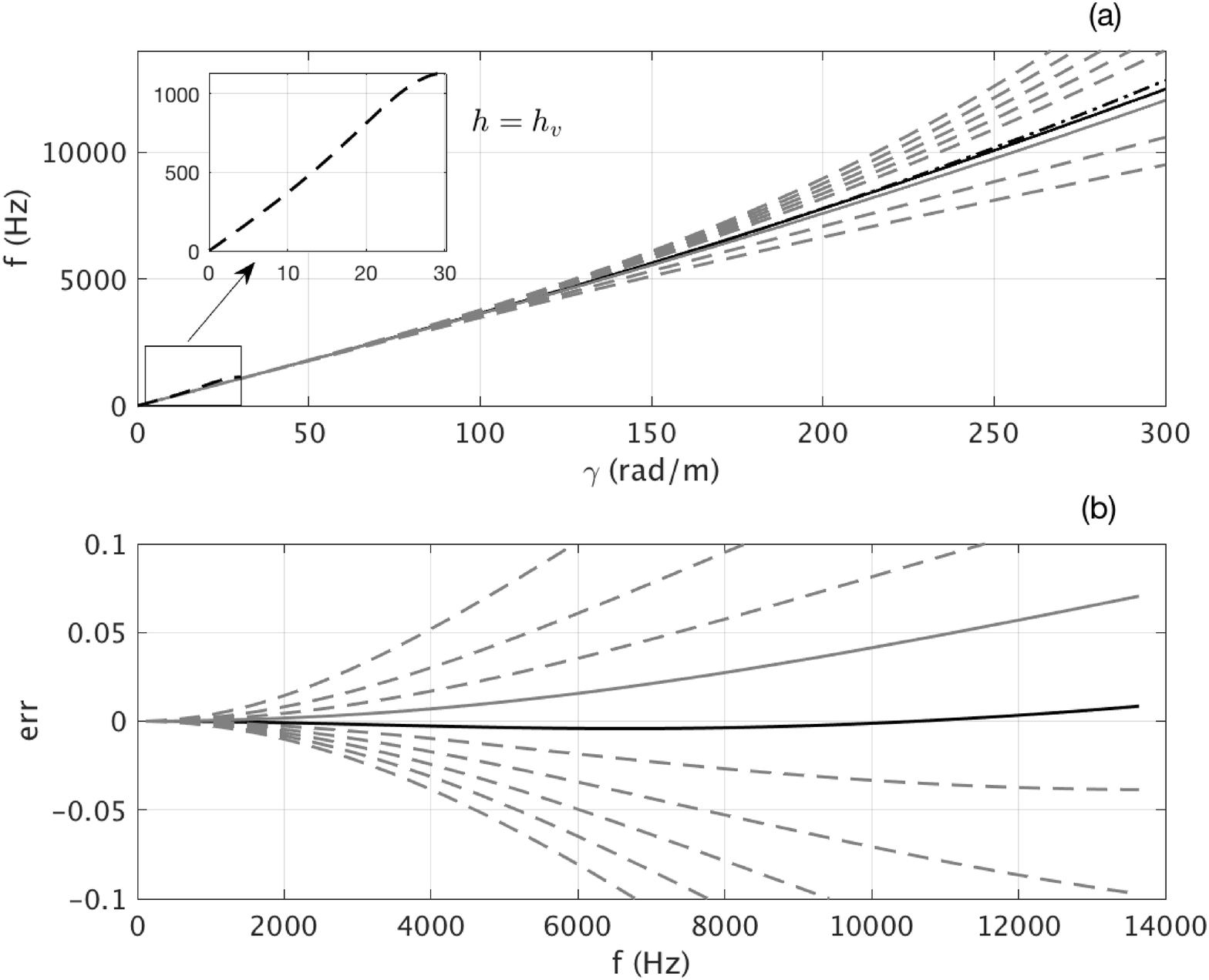}
\caption{Numerical dispersion relation and modes of the linear transverse displacement. For both figures, a string with $\rho = 8000$ kg/m$^3$, $E = 2\cdot 10^{11}$ Pa, $r = 0.2$ mm, $L = 1$ m, $T_0=50$ N was used. The sample rate is $f_s = 48$ kHz. (a): Dispersion relations. Black dash-dotted line is the continuous dispersion relation; black solid line is obtained using the stability condition \eqref{eq:StabCndStiffTheta} using $\theta_u = \bar \theta_u$ given in \eqref{eq:theta}; grey solid line is obtained with $\theta_u=1$; grey dashed lines are obtained under various other choices for $\theta_u$. The choice of the longitudinal grid spacing $h = h_v = \sqrt{E/\rho}\,\,k$ is also given in the inset. The numerical dispersion relations are obtained from  \eqref{eq:NumericalOmegs}, by transforming the difference operators in the frequency domain, as ${\bf D}^2 \rightarrow -4/h^2 \sin^2(\gamma h/2)$, ${\bf D}^4 \rightarrow 16/h^4 \sin^4(\gamma h/2)$, $\dtt \rightarrow -4/k^2 \sin^2(\omega k/2)$. Here, $\gamma$ is the wavenumber, and $\omega$ is the radian frequency (though the linear frequency $f=\omega / 2\pi$ is plotted).
(b): error on modal frequencies.
Black solid line is obtained using the stability condition \eqref{eq:StabCndStiffTheta} using $\theta_u = \bar \theta_u$ given in \eqref{eq:theta}; grey solid line is obtained with $\theta_u=1$; grey dashed lines are obtained under various other choices for $\theta_u$. In (b), the error is defined as $1-\omega/\omega^{(exact)}$. The numerical frequencies $\omega$ are as per \eqref{eq:NumericalOmegs}.}\label{fig:disprelsStiff}
\end{figure}

Consider the following modification of the linear part of \eqref{eq:FullDisc1}
\begin{equation}\label{eq:StiffFDtheta}
\left(\rho A {\bf R}(\theta_u)\dtt - T_0 {\bf D}^{2} + EI {\bf D}^{4}   \right) {\bf u}^n = 0,
\end{equation}
where the parameterised operator ${\bf R}\left(\theta_{u}\right)$ has the form
\begin{equation}\label{eq:Gdef}
{\bf R}(\theta_u) = {\bf I} + \frac{(1-\theta_u)h^2}{2}{\bf D}^{2}.
\end{equation}
One can select $\theta_u$ so that the number of grid intervals matches the number of eigenfunctions whose eigenfrequency is found below the Nyquist limit $1/2k$, resulting in  wideband dispersion reduction \cite{bilbao_book} (Chapter 7). Let such number of eigenfunctions be $N_u$. Under simply-supported conditions, 
\begin{equation}
N_u = \frac{L}{{\pi}}\sqrt{\frac{-T_0+ \sqrt{T_0^2 + \frac{4 \pi^2}{k^2} \rho E A I}}{2 EI}}
\end{equation}
Then, one selects (see \cite{ducceschi_WaveMotion_2019})
\begin{equation}\label{eq:theta}
\theta_u = \bar \theta_u \triangleq \frac{1}{2}+\frac{T_0k^2\bar h^2 + 4EI k^2}{2 \rho A \bar h^4}, \quad \text{with } \bar h = \frac{L}{N_u}.
\end{equation}
The numerical transverse eigenfrequencies $\omega_{m}$, $1 \leq m \leq N-1$ of the parameterised scheme can be determined through use of the ansatz ${\bf u}^{n} = e^{j\omega n k}\hat{{\bf u}}$, for a constant vector $\hat{{\bf u}}$, and angular frequency $\omega$, as
\begin{equation}\label{eq:NumericalOmegs}
\{\omega_m \}= \frac{2}{k}\arcsin\left(\frac{k}{2 }\text{eig}\left[{\bf R}(\theta_u)^{-1}\left( -\frac{T_0}{\rho A} {\bf D}^{2}  + \frac{EI}{\rho A} {\bf D}^{4}  \right)\right]^{1/2}\right).
 \end{equation}
 Fig. \ref{fig:disprelsStiff} shows the numerical dispersion relations and the modal frequencies for various choices of the parameter $\theta_u$. It can be seen that, under the choice \eqref{eq:theta}, the warping effects are minimised across the spectrum, though order-accuracy is not formally decreased (since the correction factor is itself $O(h^2)$), when the stability condition is fulfilled close to the stability limit. The same figure shows the results obtained when one chooses $h$ as per the natural longitudinal grid spacing $h_v = \sqrt{E/ \rho}\,\, k$:  frequency warping effects are evident, and one is not able to resolve frequencies above 1 kHz, using reference audio rates. This highlights the benefits of the choice of the nonlinear potential function $\phi$ as per \eqref{eq:PhiContDef}, which leaves the linear part of the transverse waves unaffected. Convergence of the eigenfrequencies \eqref{eq:NumericalOmegs}, as well as space-time convergence curves for the parameterised scheme \eqref{eq:StiffFDtheta} are given in \ref{app:thetaErrs}.
 
\begin{figure}
\centering
\includegraphics[width=0.92\linewidth,clip, trim = {0cm 0cm 1cm 2cm}]{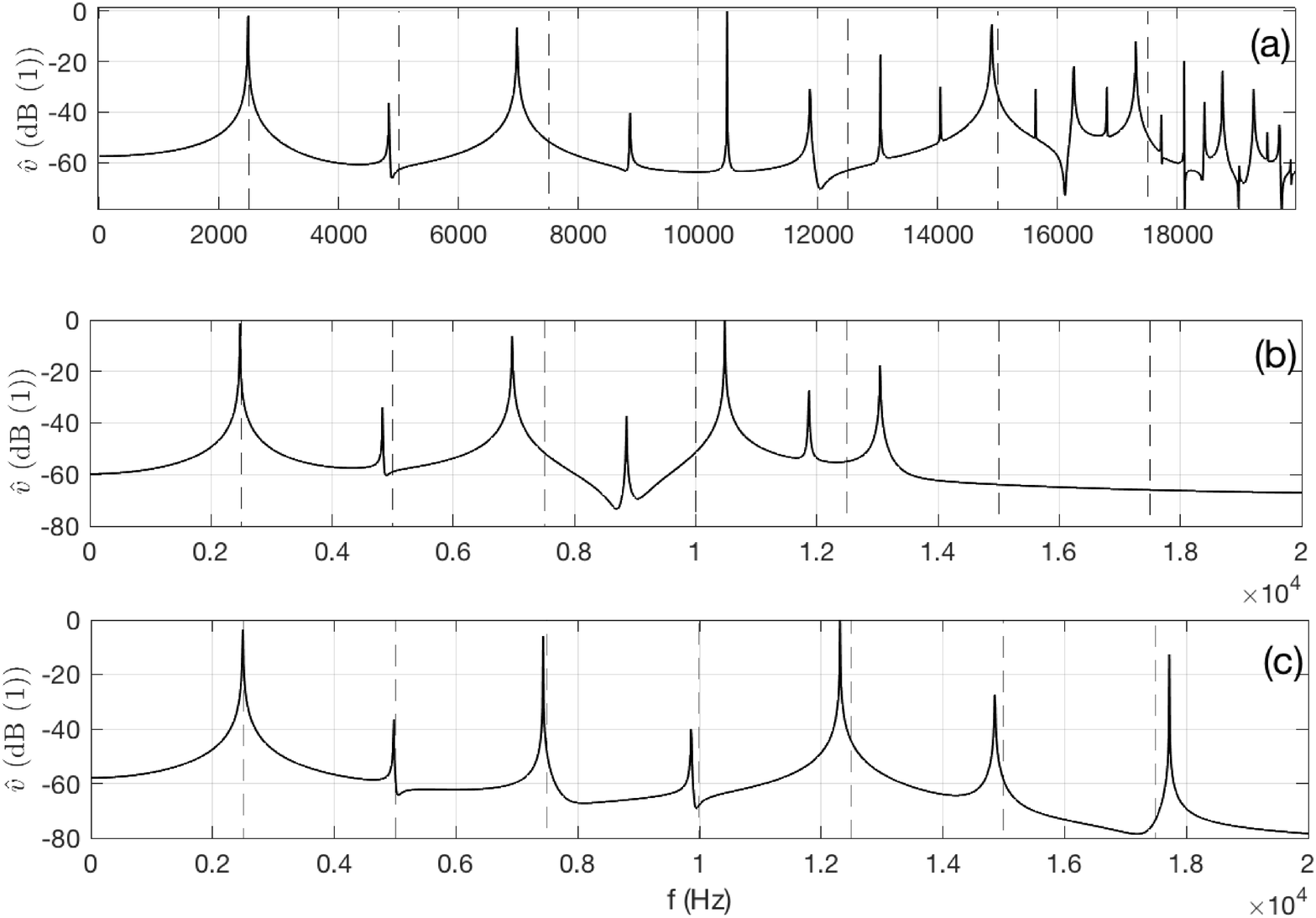}
\caption{Output spectra of the linear longitudinal displacement. For all figures, a string with $\rho = 8000$ kg/m$^3$, $E = 2\cdot 10^{11}$ Pa, $r = 0.2$ mm, $L = 1$ m, $T_0=50$ N was used. The sample rate is $f_s = 48$ kHz. (a). $N_s = (2L/\pi k)\sqrt{\rho A / T_0}$, $\theta_v = 1$. (b) $N_s = (2L/\pi k)\sqrt{\rho / E}$, $\theta_v = 1$. (c). $N_s = (2L/\pi k)\sqrt{\rho / E}$, $\theta_v$ as per \eqref{eq:thetav}. Vertical dashed lines are eigenfrequencies of the continuous system. In both cases, the longitudinal displacement was initialised with a Gaussian distribution of the kind $v(0,x) = V_0 e^{-\frac{(x/L - 1/2)^2}{2 \sigma^2}}$ with $V_0 = 1$ mm and $\sigma = 0.2$. }\label{fig:modesLong}
\end{figure}

Similarly, one may modify the linear part of \eqref{eq:FullDisc2}, expressed by means of the modal coordinates $\bf s$, in the following way
\begin{equation}\label{eq:SpectralScheme}
\left(\rho A {\bf S}(\theta_v)\dtt + T_0 {\bf \Lambda}   \right) {\bf s}^n = 0,
\end{equation}
Here,
\begin{equation}\label{eq:Sdef}
{\bf S}(\theta_v)  = {\bf I}  - \frac{(1-\theta_v)k^2}{2}{\bf \Lambda}.
\end{equation}
In this case, the numerical eigenfrequencies are given by
\begin{equation}
\{\omega_m \}= \frac{2}{k}\arcsin\left(\frac{k}{2 }\text{eig}\left[{\bf S}(\theta_v)^{-1}\left( \frac{T_0}{\rho A} {\bf \Lambda}   \right)\right]^{1/2}\right).
 \end{equation}
One can see, from Fig. \ref{fig:modesLong}, that the modal frequencies can be computed almost exactly under a suitable choice of the free parameter $\theta_v$. The value employed in the example of Fig. \ref{fig:modesLong} is
\begin{equation}\label{eq:thetav}
\theta_v = 1 + \frac{2(T_0-EA)}{7 \rho A},
\end{equation}
and it was obtained after a numerical search for the minimum of the error: $\text{max}(|\omega_m - \omega_m^{(exact)}|)$, $1\leq m \leq N_s$. Panel (a) of Fig. \ref{fig:modesLong} shows the severe modal warping effects obtained when one chooses a bound on the number of modes matching the transverse wave velocity. This example highlights the benefits of the current setup: one is indeed able to choose appropriate bounds for the finite difference grid spacing and for the number of longitudinal modes, so that frequency warping effects are minimised in both directions of motion.

The discrete energy is modified by the introduction of the operators ${\bf R}(\theta_u), {\bf S}(\theta_u)$. Corrections of the order $O(h^2), O(k^2)$ appear in the kinetic term \eqref{eq:DiscEnergy1}, which now reads
\begin{equation*}
\mathfrak{E}_k = \frac{\rho A}{2}\left( \norm{\dtm {\bf u}}^2 + \norm{\dtm {\bf v}}^2 + \frac{(\theta_u-1)h^2}{2} \norm{{\bf D}^-\dtm{\bf u}}^2  + \frac{(\theta_v-1)k^2}{2} \norm{{\bf D}^-\dtm{\bf v}}^2 \right).
\end{equation*}
This imposes a modification of the stability condition \eqref{eq:StabCndStiff}, which now is
\begin{equation}\label{eq:StabCndStiffTheta}
h \geq  h_0(\theta_u) \triangleq \sqrt{\frac{  T_0k^2  + \sqrt{ ( T_0k^2 )^2 + 16 (2\theta_u -1 )\rho A E I k^2 } }{2\rho A (2\theta_u -1)}}.
\end{equation}
Furthermore, the number of longitudinal modes should be chosen such that
\begin{equation}\label{eq:NsTheta}
N_s \leq \frac{2L}{\pi k}\sqrt{\frac{\rho A}{2(1-\theta_v)\rho A+ T_0}}.
\end{equation}
In the above, it is implied that $\theta_u > 1/2$, $2(1-\theta_v)\rho A+ T_0 > 0$.

Though the values for $\theta_u$, $\theta_v$ were  selected in the lossless case, losses are small for musical strings, and hence the wideband dispersion properties are mostly unaffected by it.
\clearpage 

\section{Application: the nonlinear struck, damped string}\label{sec:FullString}
The stiff string model detailed in Section \ref{sec:NlinWE} is here extended to include simple loss and source terms, for illustrative purposes. The equations of motion of such system, an extension of  \eqref{eq:VectorStringPsi}, are
\begin{subequations}\label{eq:StruckCnt}
\begin{align}
\left(\rho A\partial_t^2  - T_0 \partial_x^2 + EI \partial_x^4 + 2\rho A \sigma^0_u \partial_t - 2\rho A \sigma^1_u \partial_t\partial_x^2 \right) u(t,x) &= \partial_x\left( g_u \psi \right) + \mathcal{J}f(t), \label{eq:StruckCnt1} \\
\left(\rho A\partial_t^2  - T_0 \partial_x^2  + 2\rho A \sigma^0_v \partial_t  \right) v(t,x) &= \partial_x\left( g_v \psi \right). \label{eq:StruckCnt2} 
\end{align}
\end{subequations}
\eqref{eq:VectorStringPsi3} and \eqref{eq:gugv} hold here as well.
This simplified loss model depends on three parameters $\sigma_u^0,\sigma_v^0,\sigma_u^1$, the latter of which controls a frequency dependent decay. The left-hand side of \eqref{eq:StruckCnt1} is a model of the linear stiff string with loss, used in many previous works, see e.g. \cite{bilbao_book} (Chapter 7). 
Here, $\mathcal{J}$ describes the spatial distribution of the source term. For simplicity, one can choose
\begin{equation}
\mathcal{J} = \delta(x-x_{f}),
\end{equation}
where $\delta$ represents a Dirac delta function, and where $0 < x_{f} < L$ is the point of contact along the string where the source acts (not including the end points). The function $f = f(t)$ is the time evolution of the source term. A simple model for striking/plucking may be obtained by means of raised/half raised cosine \cite{bilbao_book}, as
\begin{equation}
f(t) = \frac{F_s}{2}\left( 1 - \cos (\zeta \pi (t-t_0))/t_s\right), \,\,\, t_0 \leq t \leq t_0+t_s,
\end{equation}
and $f(t) = 0$ elsewhere. The parameters $F_s$ (in Newtons), $t_s$ (in seconds) and $\zeta \in \{ 1,2 \}$ control, respectively, the maximum exerted force, the contact duration, and the source type (1 for pluck, 2 for strike). 

System \eqref{eq:StruckCnt} maintains a notion of passivity, in that it is possible to obtain an energy balance of the kind
\begin{equation}\label{eq:dHdtCnt}
\frac{dH}{dt} = -2 \rho A \mathcal{P}(t) + \partial_t u(t,x_f)f(t),
\end{equation}
where the energy $H(t)$ has the form \eqref{eq:HamWEPsi}, and where
\begin{equation}
\mathcal{P}(t) = \sigma_u^0 \norm{\partial_t u}^2 + \sigma_v^0 \norm{\partial_t v}^2 + \sigma_u^1 \norm{ \partial_t\partial_{x} u}^2 \geq 0.
\end{equation}
The time evolution of the energy \eqref{eq:dHdtCnt} expresses the power balance of the string with loss and source. In particular, there is no autonomous production of energy within the system: power is dissipated according to $2\rho A\mathcal{P}(t)$, and injected according to $\partial_t u(t,x_f)f(t)$.

\subsection{Discretisation}

The Dirac delta function can be approximated here by a vector ${\bf J} = J_m$, $1 \leq m \leq N-1$, i.e. \cite{bilbao_book} (Chapter 5)
\begin{equation}
\mathcal{J} = \delta(x-x_f) \rightarrow {\bf J},
\end{equation}
where 
\begin{subequations}
\begin{align}
    J_{m_f} &= (1-\alpha)/h, \,\, J_{m_f+1} = \alpha / h,
\end{align}
\end{subequations}
Here $\alpha = x_f/h -  m_f$, $m_f = \text{floor}(x_f/h)$. ${\bf J}$ is zero elsewhere.  The loss terms are discretised simply as
\begin{equation}
\sigma^0_u \partial_t - \sigma^1_u \partial_t\partial_x^2  \rightarrow  \sigma^0_u \dtd - \sigma^1_u \dtd {\bf D}^2, \quad \sigma^0_v \partial_t \rightarrow \sigma_v^0 \dtd .
\end{equation}
The source is simply approximated as $f(t) \rightarrow f^n$, though one should be aware that formal order of accuracy may not be preserved under such choice \cite{chabassier_CRM_2017}. This particular discretisation leads to passivity in the numerical setting. Via the same steps leading to \eqref{eq:FullDiscEner} in the conservative case, one gets
\begin{equation}\label{eq:dHdtFD}
\dtp {\mathfrak{h}^{n-1/2}} = -2\rho A \mathfrak{p}^{n} + \innp{{\bf J},\dtd{\bf u}^{n}} f^n,
\end{equation}
where
\begin{equation}
\mathfrak{p}= \sigma_u^0\norm{\dtd {\bf u}}^2 + \sigma_v^0\norm{\dtd {\bf v}}^2 + \sigma_u^1\norm{ {\bf D}^{-} \dtd{\bf u}}^2 \geq 0.
\end{equation}
Notice that stability condition \eqref{eq:StabCndStiffTheta} must be enforced here too.

\subsection{Experiments}
The performance of the discrete system is now checked in a number of numerical experiments. For those, the string and source parameters are selected as in Table \ref{tab:strParams}.

\begin{table}
\begin{center}
\begin{tabular}{ m{3cm} m{2cm}  m{3cm} }
& units  & value  \\
  \cmidrule(l){1-3} 
$\rho$ 		&kg/m$^3$			&8000              \\
$E$ 		&Pa 				&$2\cdot 10^{11}$  \\
$T_0$		&N 			 		&40                \\
$L$			&m 			        &1	               \\
$r$ 		&mm 			    &$0.29$            \\
$\sigma^0_u$&1/s                &0.1               \\ 
$\sigma^0_v$&1/s                & 0.2              \\
$\sigma^1_u$&m$^2$/s            & $4\cdot 10^{-4}$ \\
            &                   &                  \\
$t_0$       &ms                 &1                 \\
$t_s$       &ms                 &0.8               \\
$x_f$       &m                  &0.72              \\
$\mu$       &-                  &2                 \\
$F_s$       &N                  &[2.5,5,7.5]       \\
\cmidrule(l){1-3} 
\end{tabular}
\end{center}
\caption{String Parameters used in the simulations.}\label{tab:strParams}
\end{table}

Notice that three possible values for $F_s$ can be selected from the table. Moreover, output is extracted at $x_o = 0.32$ m. For the simulations, a base sample rate  $f_{s0} = 48\cdot 10^3$ Hz is used, and oversampling factors are indicated in the figures. For all the simulations, the grid spacing is chosen as $h = 1.05 h_0$, where $h_0$ is the limit of stability defined in \eqref{eq:StabCndStiffTheta}. The number of longitudinal modes $N_s$ is selected as per bound \eqref{eq:NsTheta}.
In particular, in order to set $\theta_u$ in \eqref{eq:theta}, one should employ a modified value for $\bar h$, i.e.  $\bar h = L/1.05 N_u$. This is to avoid erratic behaviour very close to Nyquist, as noted in Section \ref{sec:experiments}. 
\begin{figure}
\includegraphics[width=\linewidth]{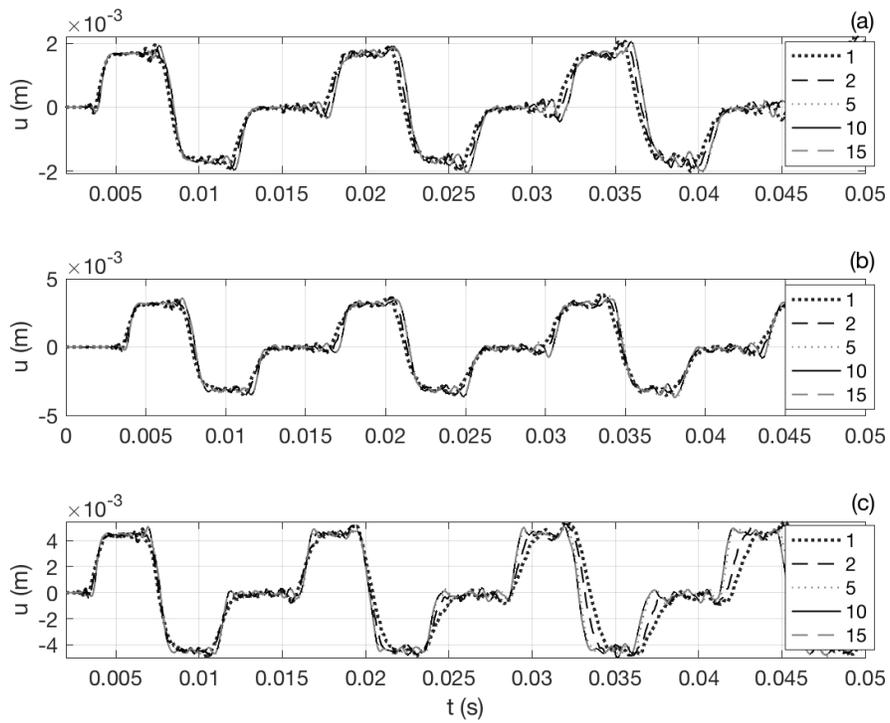}
\caption{Transverse displacement. (a): $F_s=2.5$ N; (b): $F_s=5$ N; (c): $F_s=7.5$ N. Oversampling factors as indicated.}\label{fig:TransvAll}
\end{figure}
A first experiment is presented in Fig. \ref{fig:TransvAll}. Here, the panels show the time evolution of the transverse waves for the three forcing values, for a number of oversampled solutions with varying oversampling factors. One can see that the numerical solutions tend to converge to a common solution as the sample rate increases. The waveforms show that the wave velocity (and thus the frequency) increases with the forcing amplitude, but that the overall maximum amplitude of the wavefront is approximately the same for the last two cases ($F_s=5,7.5$ N). 
\begin{figure}
\includegraphics[width=\linewidth, clip, trim={5cm 0cm 4cm 1cm}]{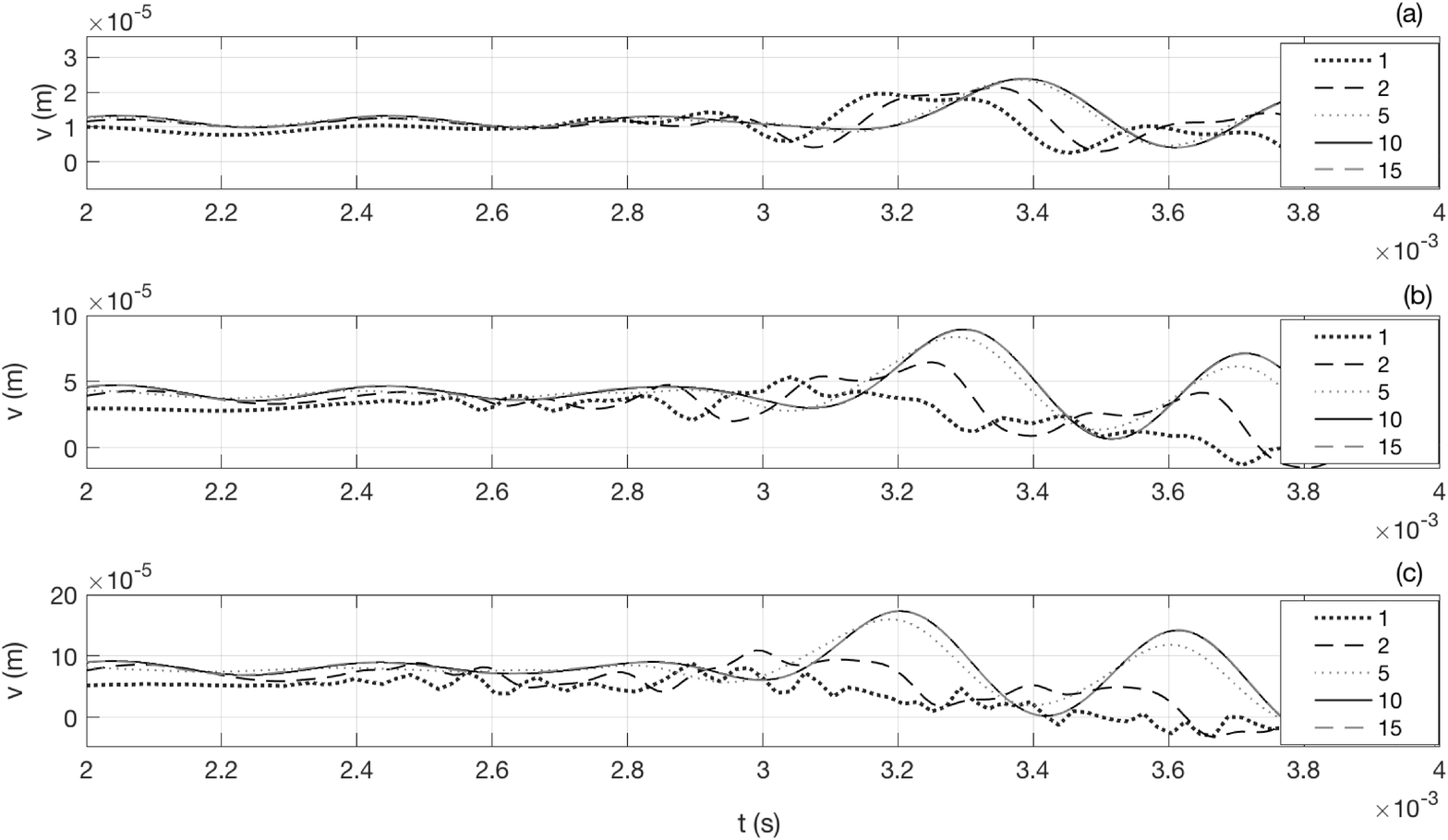}
\caption{Longitudinal displacement. (a): $F_s=2.5$ N; (b): $F_s=5$ N; (c): $F_s=7.5$ N. Oversampling factors as indicated.}\label{fig:LongAll}
\end{figure}
Input power, in fact, is here converted to longitudinal motion, rather than larger transverse amplitude, as one would observe in the linear case. These facts are evident from Fig. \ref{fig:LongAll}, where one sees a progressive increase of the amplitude of the longitudinal waves with input forcing. 
\begin{figure}[hbt]
\centering
\includegraphics[width=0.92\linewidth]{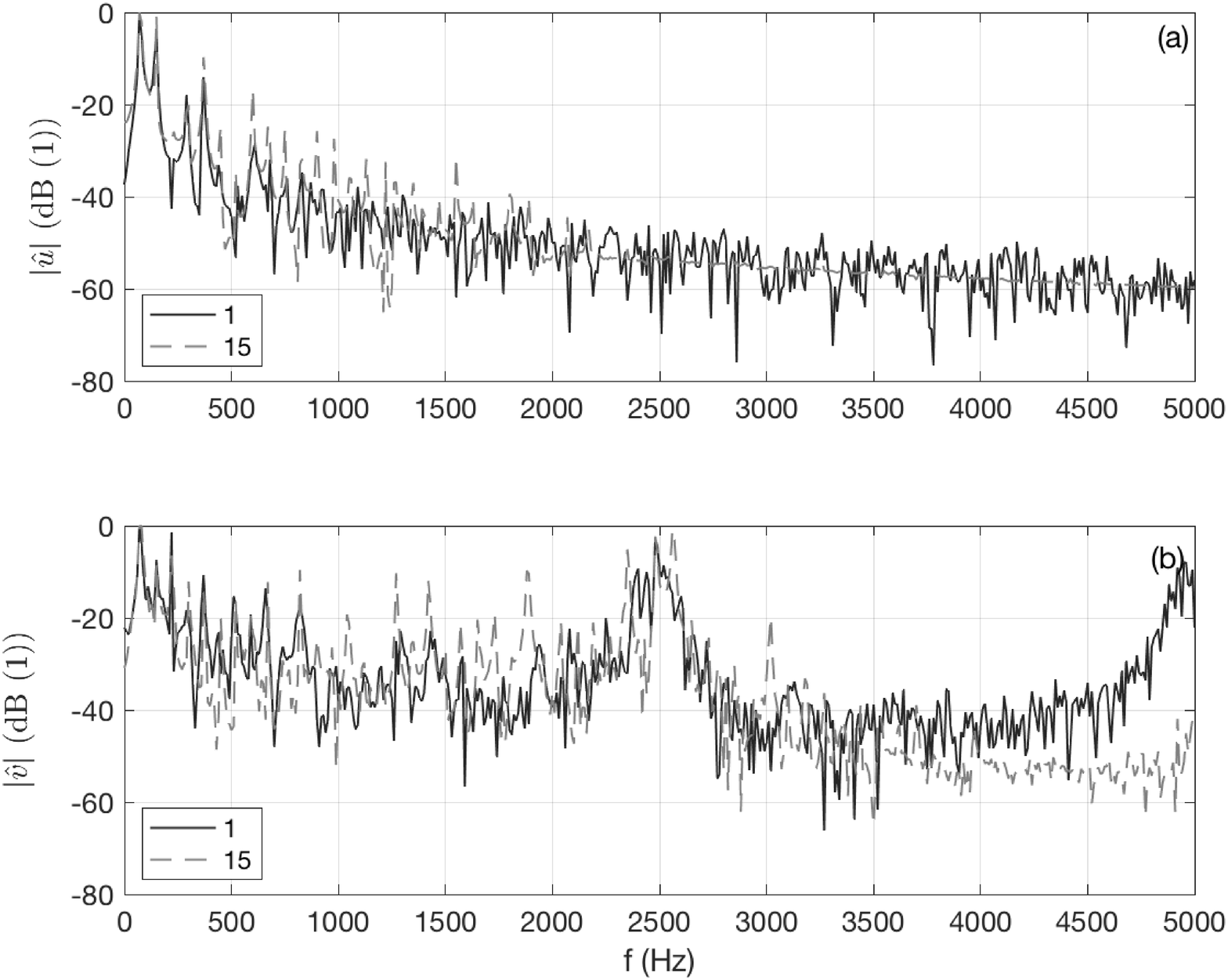}
\caption{Output spectra. Top: transverse wave spectrum. Bottom: longitudinal wave spectrum. $F_s=5$ N. Oversampling factors as indicated.}\label{fig:SpectraSource}
\end{figure}
The effects of aliasing can be observed in Fig. \ref{fig:SpectraSource}, where the spectra of the solution computed at 15 times audio rate versus the solution computed at audio rate are compared: the longitudinal spectrum shows some aliasing, though the transverse spectrum presents only some floor noise at $-50$ dB beyond 2.5 kHz.
\begin{figure}[hbt]
\centering
\includegraphics[width=0.83\linewidth,clip, trim = {2cm 3.0cm 2cm 6.0cm}]{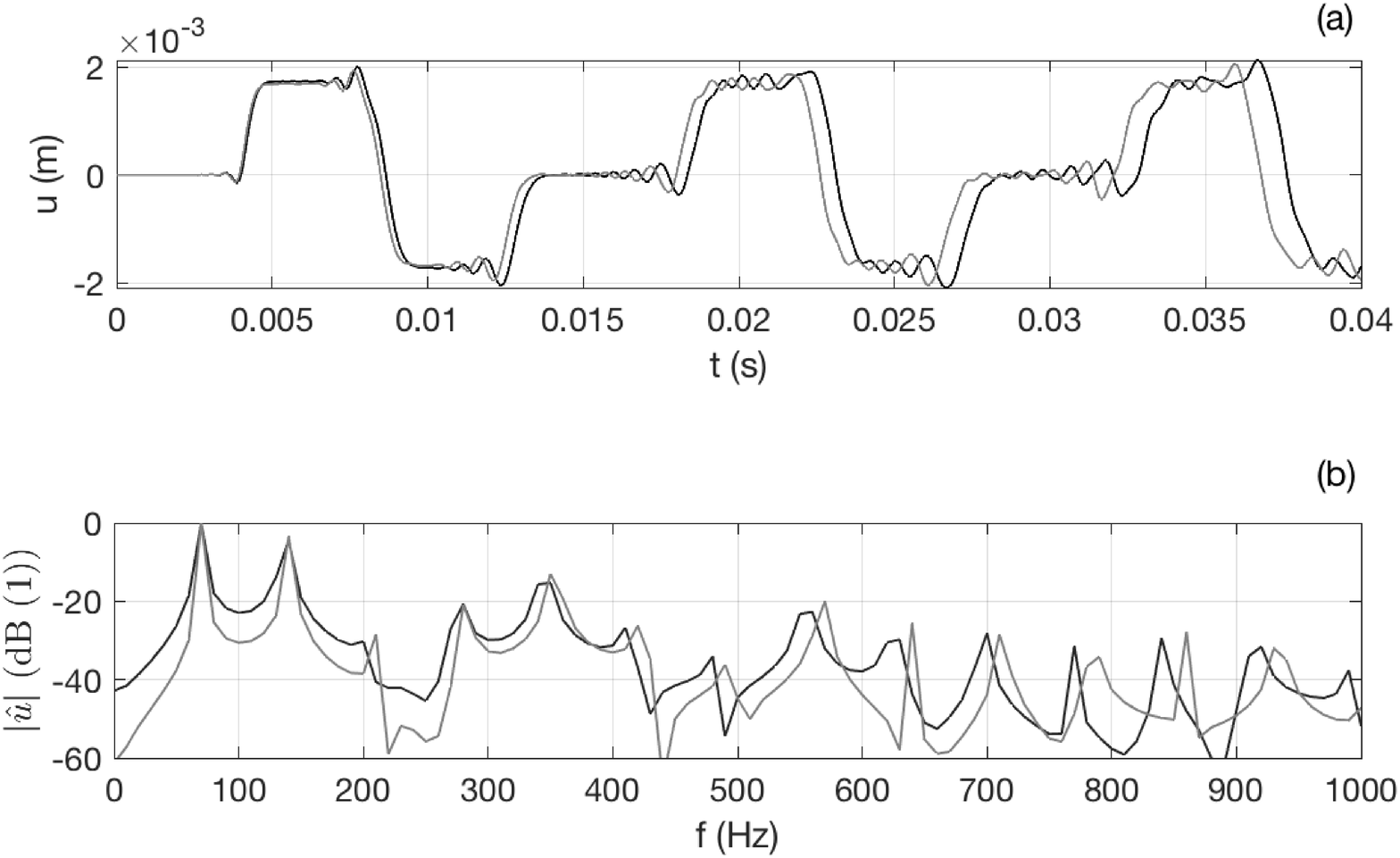}
\caption{Linear (Black) vs Nonlinear (Grey) Dynamics. Top: transverse displacement. Bottom: spectrum. $F_s=2.5$ N. Sample rate $f_s = 15f_{s0}$.}\label{fig:LinNlin}
\end{figure}
A final experiment is presented in Fig. \ref{fig:LinNlin}, where the linear vs nonlinear dynamics are compared, for the smallest forcing value. Notice that the maximum amplitude of the dynamics is small, of the order of 2 mm. One can nonetheless appreciate significant differences between the two solutions, highlighting the importance of nonlinear effects for accurate simulations of string dynamics. Sound examples and Matlab code snippets are available at the companion webpage\footnote{https://mdphys.org/nonlinearString.html}, using the parameters from Table \ref{tab:strParams}.

\section{Conclusions}
String vibration is a subject of interest in many engineering disciplines. For the purpose of sound synthesis by physical modelling, as well as for advanced musical acoustics investigation, one must include nonlinear effects. These are distributed across the string, and arise as a consequence of the string's large stretching. In this work, an application of quadratisation methods was offered. This allows to  solve for the nonlinearities using one single matrix inversion, thus bypassing the machinery of iterative methods. Longitudinal and transverse waves were solved using a mixed approach, so that the wave speeds may be resolved at reference audio rates, without the need to oversample. The use of free theta parameters allows to implement wideband numerical dispersion reduction. Stability of the numerical scheme was approached via energy conservation, and illustrated numerically. Extension of these techniques to other nonlinearities in acoustics, such as collision dynamics \cite{bilbao_AAUA_2015,ducceschi_DAFx_2019,chatziioannou_JSV_2015,ducceschi_jasa_2021}, von K{\'a}rm{\'a}n plates \cite{bilbao_book,ducceschi_PhysD_2014}, and others, is possible. Important aspects, worthy of future investigations, include a formal proof of convergence and accuracy,  an analysis of structure of the scheme's matrices for improved performance, as well as the design of schemes with reduced aliasing in the longitudinal component.

\section*{Acknowledgements}
The first author wishes to acknowledge the Leverhulme Trust, who  supported this research with a Leverhulme Early Career Fellowship. This work was also supported by the European Research Council (ERC), under the European Union's  Horizon 2020 research and innovation programme,  grant 2020-StG-950084-NEMUS. The anonymous reviewers are also thanked for their comments and suggestions.

\appendix

\section{Error curves for the quadratised scheme}\label{app:DuffErrs}

\begin{figure}
\includegraphics[width=\linewidth]{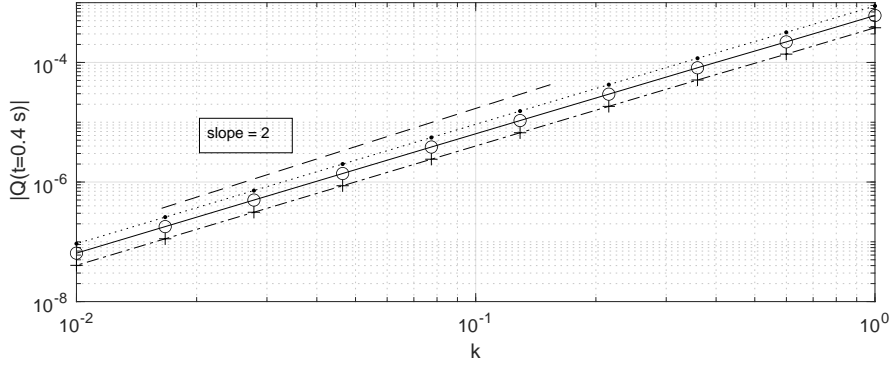}
\caption{Error curves of the quadratised scheme \eqref{eq:DuffingQuad}. Here, $+$ is $\gamma = 0.6$, $\circ$ is $\gamma = 0.8$, $\cdot$ is $\gamma = 1$. The initial displacement is selected as $u_0=3.7$, and the error is computed at the time $t=0.4$ s. The dashed line has a slope of 2.}\label{fig:DuffingSlopes}
\end{figure}

The simple case of a scalar nonlinear ODE is considered here, to check the order of accuracy of the proposed quadratised scheme. The scalar ODE reads
\begin{equation}\label{eq:Duffing}
\frac{d^2u}{dt^2} = - u - \gamma u^3.
\end{equation}
This is a type of lossless Duffing equation, with zero input. For this equation, an analytic solution exisits as
\begin{equation}\label{eq:DuffingAnalytic}
u(t) = u_0 \, \text{cn}\left( \sqrt{1 + \gamma u_0^2} \, t ; \, \frac{\gamma u_0^2}{2\gamma u_0^2 + 2}\right),
\end{equation}
where $\text{cn}(a;b)$ is the Jacobi elliptic function with argument $a$ and parameter $b$. Here, it is assumed that $u(t=0) = u_0$, $\frac{du(t=0)}{dt}=0$. Numerical integration of \eqref{eq:Duffing} using the quadratised energy proceeds as follows. First, define the auxiliary variable 
\begin{equation}
\psi = \sqrt{\frac{\gamma}{2}u^4} = \sqrt{\frac{\gamma}{2}}u^2.
\end{equation}
The scheme is constructed as
\begin{equation}\label{eq:DuffingQuad}
\dtt u^n = -u^n - g^n \mtp \psi^{n-1/2}, \,\, \dtp \psi^{n-1/2} = g^n \dtd u^n,\,\, g^n = \frac{d\psi}{du}\big|_{u=u^n}.
\end{equation}
Numerical initial conditions are specified as
\begin{equation}
u^0 = u_0, \,\, u^1 = u^0 - \frac{k^2}{2} (u^0 + \gamma (u^0)^3),  \,\, \psi^{1/2} = \sqrt{\frac{\gamma}{2}}(u^0)^2.
\end{equation}
These are second-order accurate. Note, in particular, that $u^1$ was obtained implementing $\frac{d}{dt} = \dtp - \frac{k}{2}\dtt + O(k^2)$. The error at the time $t_e = \text{round}(0.4/k)k$ is then computed as 
\begin{equation}
Q= u(t_e) - u^{t_e/k},
\end{equation}
that is, the difference between the analytic solution \eqref{eq:DuffingAnalytic}  and the output of the finite difference scheme at the corresponding time step.
The error curves, under various choices of the nonlinear parameter $\gamma$, are plotted in Fig. \ref{fig:DuffingSlopes}, and are second-order.

\section{Error curves for the linear $\theta$-scheme}\label{app:thetaErrs}

\begin{figure}[hbt]
\includegraphics[width=\linewidth]{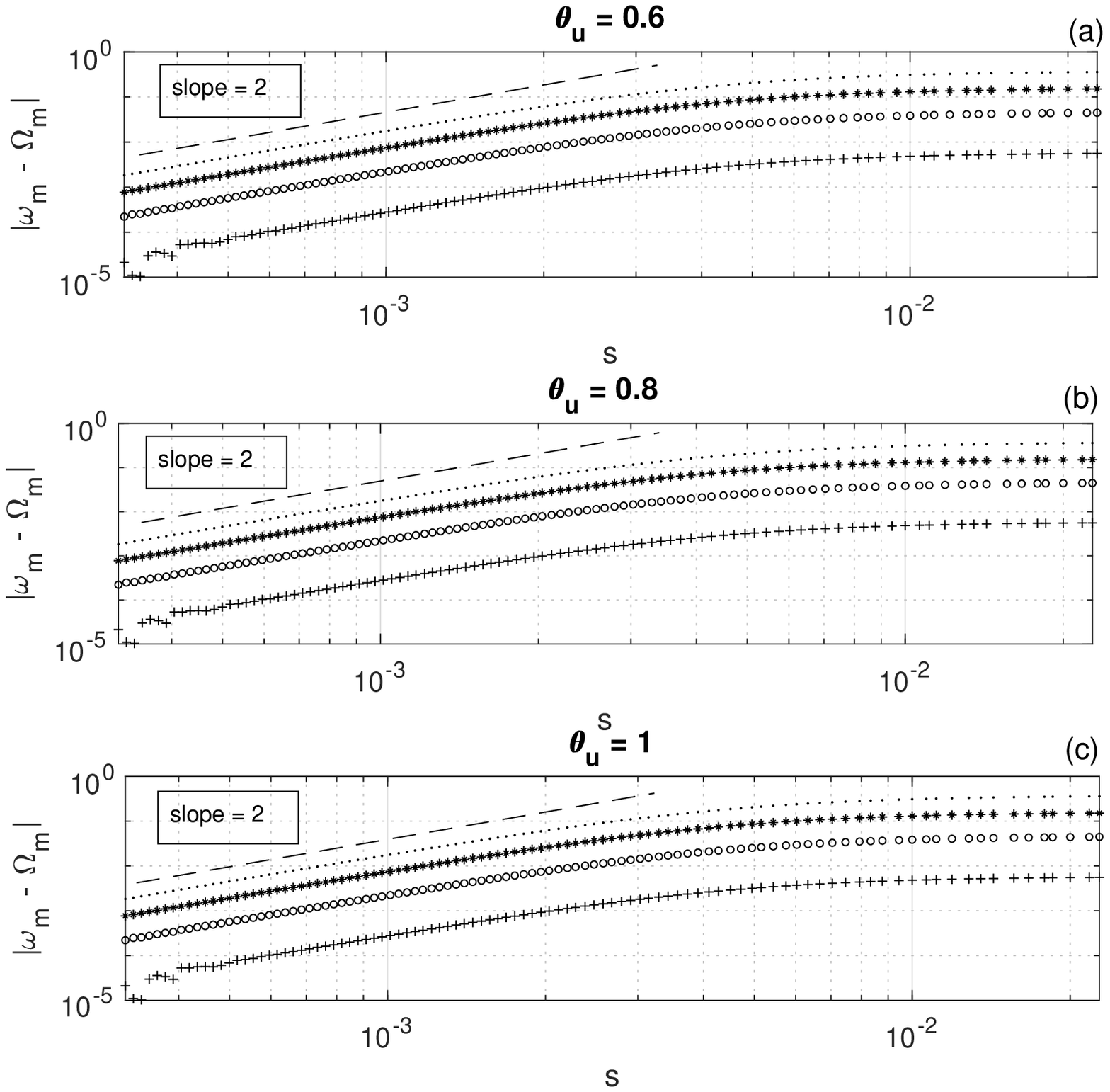}
\caption{Numerical errors of the first four eigenfrequencies of scheme \eqref{eq:StiffFDtheta}, under various choices of $\theta_u$. Here, the numerical eigenfrequencies $\omega_m$ are given by \eqref{eq:NumericalOmegs}. The analytic eigefrequencies  $\Omega_m$ are as per \eqref{eq:AnOmegs}. For all panels, $+$ is $m=1$, $\circ$ is $m=2$, $\text{*}$ is $m=3$, and $\cdot$ is $m=4$. The superimposed dashed lines have a slope of 2. A string with $\rho = 8000$ kg/m$^3$, $E = 2\cdot 10^{11}$ Pa, $r = 0.2$ mm, $L = 1$ m, $T_0=50$ N was used for the simulations.}\label{fig:slopesEigen}
\end{figure}

The order of accuracy of the parameterised scheme \eqref{eq:StiffFDtheta} is checked. From the stability condition \eqref{eq:StabCndStiffTheta}, it is clear that the grid spacing $h$ is not simply proportional to the time step $k$, as would be in a CFL-like condition. For this reason, the specification of the order of accuracy may be done in terms a new variable that parametrises uniquely the level curve $h-h_0(k,\theta_u) \approx 0$ in the $(h,k)$ plane, where $h_0$ is as per \eqref{eq:StabCndStiffTheta}. First, note that level curve can be expressed as
\begin{equation}\label{eq:TmStp}
k = h^2  \sqrt{\frac{\rho A(2\theta_u-1)}{T_0h^2+4 EI}}.
\end{equation}
This curve is naturally parametrised by the arclength $s$, defined as
\begin{equation}
s = \int_0^h \sqrt{1 + \left( \frac{dk}{dh}\right)^2}\,\, \mathrm{d}h = h + \frac{\rho A(2\theta_u-1)}{6 EI} h^3 + O(h^5).
\end{equation}
In Fig. \ref{fig:slopesEigen}, a check on the convergence of the numerical eigenfrequencies is performed. It is seen that clearly the errors are $O(s^2)$, for small enough $s$, and irrespective of the choice of $\theta_u$ (so long, of course, that $\theta_u>1/2$, as specified at the end of Section \ref{subsec:NumDisp}).  

A second check is performed in the time domain. Here, the value of the numerical solution (at a given location along the grid, and at a given time) is compared to an analytic solution. The analytic solution is constructed as follows. Consider a modal expansion of the solution of the linear part of \eqref{eq:VectorString1a} (i.e., neglecting the right-hand side). This is
\begin{equation}
u(t,x) = \sum_{m=1}^M \left( C_m e^{j\Omega_m t} + C^\star_m e^{-j\Omega_m t}\right) \sin \frac{m \pi x}{L}.
\end{equation}
Here, $C_m \in \mathbb{C}$ is a complex constant, with complex conjugate $C_m^\star$, to be determined from the initial conditions. $M$ is, in theory, infinite, but for all practical puroposes it must be truncated to an integer (specified below). 
The resonant frequencies are given by 
\begin{equation}\label{eq:AnOmegs}
\Omega_m = \sqrt{\frac{T_0}{\rho A} \left(\frac{m \pi}{L} \right)^2 + \frac{EI}{\rho A} \left(\frac{m \pi}{L} \right)^4}.
\end{equation}
The initial displacement is selected as a centered raised cosine with compact support:
\begin{equation}
u(t=0) \triangleq u_0(x)  = 
\left\{ 
\begin{array}{ll}
1 - \cos \left(\frac{2\pi(x-L/4)}{L/2}\right)& \text{if }L/4\leq x \leq 3L/4, \\
0 & \text{elsewhere. }
\end{array}
\right.
\end{equation}
The initial velocity is selected as $\frac{du(t=0)}{dt}=0$. Under such initial conditions, the solution is
\begin{equation}\label{eq:ddd}
u(t,x) = \sum_{m=1}^{M}2C_{m}\cos(\Omega_m t) \sin \frac{m \pi x}{L},
\end{equation}
with
\begin{equation}
C_{m}  = \frac{\left< u_0, \sin\frac{m \pi x}{L} \right>}{L},
\end{equation}
where the inner product is defined in \eqref{eq:innpCnt}. These integrals can be computed analytically, ultimately yielding an analytic expression for $u(t,x)$ for all times, when substituted back into \eqref{eq:ddd}.

\begin{figure}[hbt]
\includegraphics[width=\linewidth]{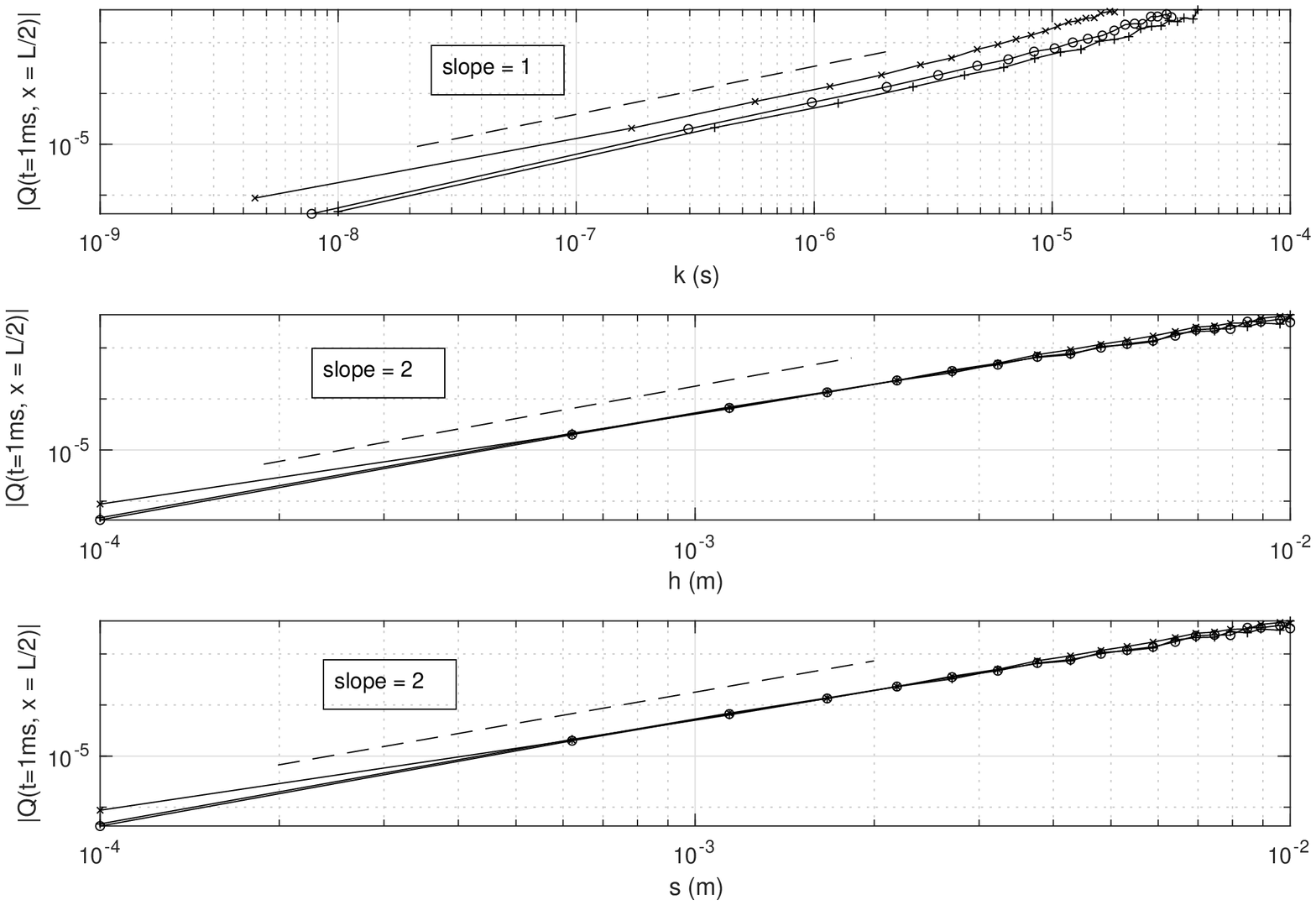}
\caption{Space-time error curves for scheme \eqref{eq:StiffFDtheta}, under various choices of the parameter $\theta_u$, and plotted against $k$, $h$ and $s$. For all panels $x$ is $\theta_u=0.6$, $\circ$ is $\theta_u = 0.8$, and $+$ is $\theta_u=1$. Dashed lines give the error slopes. A string with $\rho = 8000$ kg/m$^3$, $E = 2\cdot 10^{11}$ Pa, $r = 0.2$ mm, $L = 1$ m, $T_0=50$ N was used for the simulations.}\label{eq:SpaceTimeErrs}
\end{figure}

The numerical scheme is run as follows. First, the numerical initial conditions are implemented as
\begin{subequations}
\begin{align}
{ u}^0_m &= u_0(x=mh), \\
{\bf u}^1 &= {\bf u}^0 + \frac{k^2}{2} \left( \frac{T_0}{\rho A}  {\bf D}^2 - \frac{EI}{\rho A} {\bf D}^4\right) {\bf u}^0. \label{eq:IC2}
\end{align}
\end{subequations}
These expressions for the initial conditions are second-order accurate. In particular, \eqref{eq:IC2} is obtained using the fact that $\frac{d}{dt} = \delta_{t+} - \frac{k}{2}\delta_{tt} +O(k^2)$. Then, a number $N$ of subintervals is selected, with $N$ even. The corresponding grid spacing $h$ is obtained as $h = L/N$, and the time step is then selected as \eqref{eq:TmStp}. The error is computed at $x_e = L/2$, $t_e = \text{round}(10^{-3}/k)k$ (i.e. 1 millisecond), as 
\begin{equation}\label{eq:ErrSpaceTime}
Q = u(t_e,x_e) - u^{t_e/k}_{N/2 + 1},
\end{equation}
that is, the difference between the analytic solution \eqref{eq:ddd}, and the output of the finite difference scheme at the corresponding grid point and time step. The analytic solution is computed using a total number of modes $M = 10^6$. The results are plotted in Fig. \ref{eq:SpaceTimeErrs}. Note that the error is $O(s^2)$ (as expected). Since the leading term in the expansion for $s$ is equal to $h$, plotting the error against the grid spacing gives again a slope of 2. On the other hand, since $k \propto h^2 $ for small $h$, the error has a slope of $1$ when plotted against the time step. Note that these results are consistent with the error curves presented in \cite{ducceschi_WaveMotion_2019}.

\section{Comparing implicit and explicit schemes for the linear wave equation with stiffness}\label{app:ExpVsImp}

The local truncation error of the explicit scheme \eqref{eq:StiffFDtheta} (with $\theta_u=1$), is $O(k) + O(h^2)$, and is thus  first-order in time. Here, a comparison against an implicit scheme, formally second-order in space and time, is given. To that end, consider the following implicit scheme
\begin{equation}\label{eq:ImpScheme}
    \left(\rho A {\bf I} + \frac{k^2 EI}{2}{\bf D}^4\right) \delta_{tt} {\bf u}^n = T_0 {\bf D}^2 {\bf u}^n - EI  {\bf D}^4{\bf u},
\end{equation}
which conserves the discrete energy 
\begin{equation}
{\mathfrak h}^{n-1/2} = {\mathfrak E}_k^{n-1/2} + {\mathfrak E}_{pl}^{n-1/2},
\end{equation}
where
\begin{subequations}
  \begin{align}
    {\mathfrak E}_k^{n-1/2} &= \frac{\rho A}{2}\|\delta_{t-}{\bf u}^n \|^2 + \frac{EI k^2}{4}\| {\bf D}^2\delta_{t-}{\bf u}^n\|^2 , \\
    {\mathfrak E}_{pl}^{n-1/2} &= \frac{T_0}{2}\left<{\bf D}^-{\bf u}^n,{\bf D}^-\etm{\bf u}^{n}\right> + \frac{EI}{2}\left<{\bf D}^2{\bf u}^n,{\bf D}^2\etm{\bf u}^{n}\right> .
    \end{align}
\end{subequations}
The energy is non-negative (and therefore stability ensues) if and only if
\begin{equation}
    h \geq h_{imp} \triangleq \sqrt{T_0/\rho A}\, k,
\end{equation}
that is, the CFL condition. The local truncation error is $\tau = O(k^2) + O(h^2)$, that is, second-order in $k$ and $h$. This scheme is now compared to the explicit scheme \eqref{eq:StiffFDtheta} (using $\theta_u=1$), that is
\begin{equation}\label{eq:ExpScheme}
    \rho A \delta_{tt} {\bf u}^n = T_0 {\bf D}^2 {\bf u}^n - EI  {\bf D}^4{\bf u},
\end{equation} for which a stability condition arises as per \eqref{eq:StabCndStiff}, i.e. 
\begin{equation}
    h\geq h_{0}=\sqrt{\frac{  T_0k^2  + \sqrt{ ( T_0k^2 )^2 + 16 \rho A E I k^2 } }{2\rho A }}.
\end{equation}
Since $h_{0} > h_{imp}$, scheme \eqref{eq:ExpScheme} is more efficient than \eqref{eq:ImpScheme}, when the same sample rate is used. Efficiency is also greatly improved by the fact that  \eqref{eq:ExpScheme} is explicit, whereas \eqref{eq:ImpScheme} requires the solution of a pentadiagonal system at each time step. The question that is addressed here is whether formal first-order accuracy in time has in fact a negative impact on the convergence curves of the proposed explicit scheme. A number of checks are now performed to address this question. The checks will be performed either selecting the sample rate at the input and computing the grid spacing at the limit of stability, or selecting the grid spacing at the input and using the corresponding sample rate at the limit of stability. It will be clear that the explicit scheme outperforms the implicit one.
\begin{figure}[hbt]
    \centering
    \includegraphics[width=\linewidth]{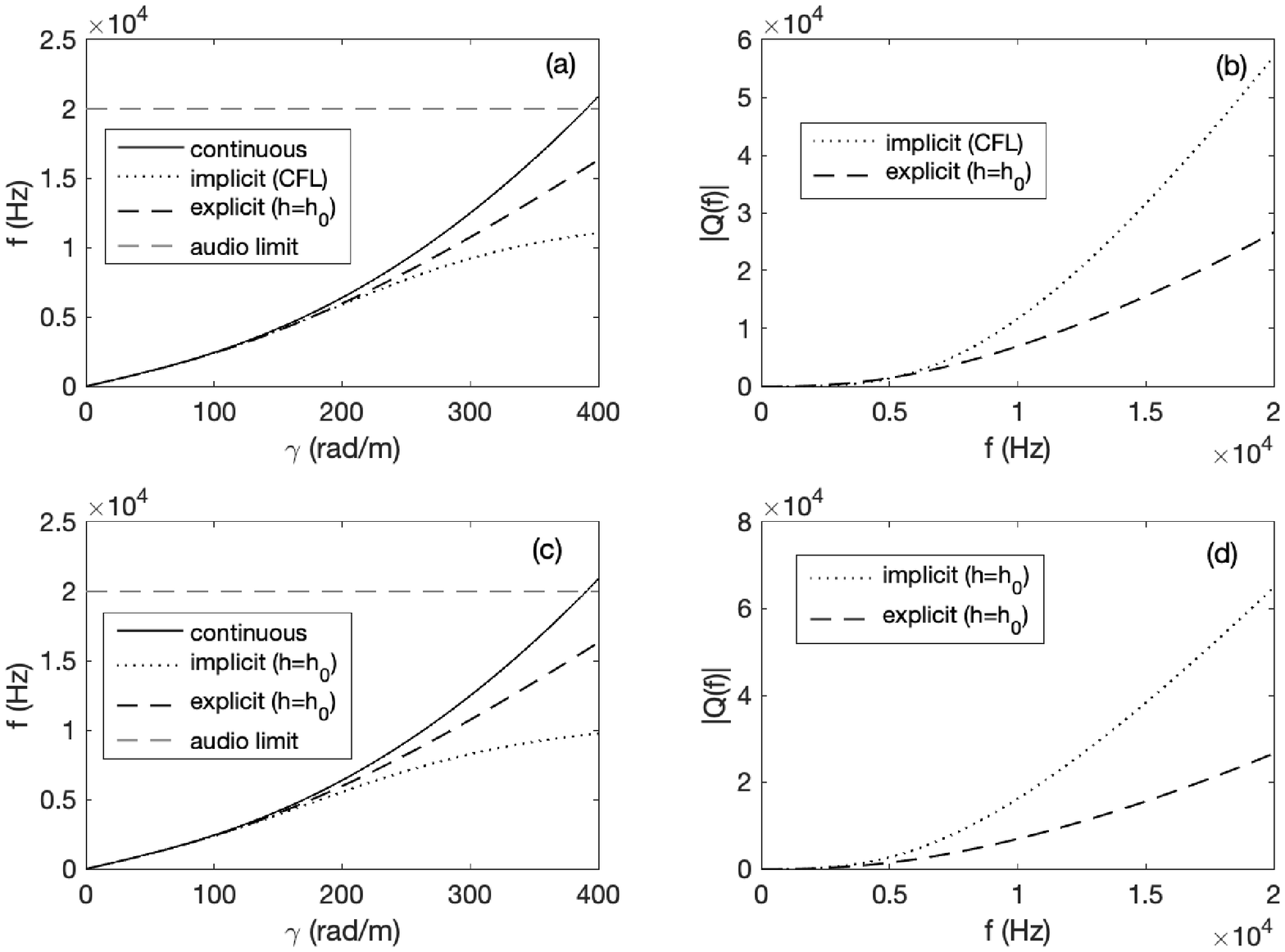}
    \caption{Continuous and numerical dispersion relations. Panels (a) and (b) plot the dispesion relation and the error $E$ (i.e. the difference between the continuous and numerical dispersion relation) for the schemes \eqref{eq:ImpScheme} and \eqref{eq:ExpScheme}, when they are both run at their respective limits of stability. The parameters are $f_s = 48 \cdot 10^3$, $\rho = 8000$ kg/m$^3$, radius $r=0.3$ mm, $T_0 =40$ N and $E = 200$ GPa. Panels (c) and (d) repeat the same experiments, except that the implicit scheme is run at the limit of stability of the explicit scheme. Note that, in both cases, the explicit scheme has lower numerical dispersion. It is also much more efficient, since in the first case, the grid points of the implicit scheme are more than twice as many, and the solution of a pentadiagonal linear system is required at each time step. }\label{fig:Dispersion}
\end{figure}

First, consider Fig. \ref{fig:Dispersion}: there, the numerical dispersion relations are plotted for the two schemes, and compared to the continuous dispersion relations. The sample rate is selected as input parameter. Two experiments are presented. In the first experiment (panels (a) and (b)), the two schemes are run at their respective limits of stability, i.e. $h=h_{imp}$ for \eqref{eq:ImpScheme} and $h=h_{0}$ for \eqref{eq:ExpScheme}. It is seen that the explicit scheme has much lower numerical dispersion than the implicit scheme, particularly for higher frequencies, {despite using a number of grid points that is less than half compared to the implicit scheme.} Here, the number of eigenfrequencies of the continuous system below 20 kHz is 139; the number of grid points of the explicit scheme is 168, and the number of grid points of the implicit scheme is 360. The second experiment (panels (c) and (d)) shows a slightly fairer comparison, using the \emph{same} sample rate and grid spacing for both schemes. Here, the performance of the implicit scheme is in fact worse than in the first experiment, as the error in the frequency domain becomes bigger. 
The implicit discretisation, when run at the limit of stability, is wasteful, since the extra grid points correspond to numerical modes that are completely warped. Matlab simulations in the time domain show that runtime of the implicit scheme is roughly 10 times longer than the explicit one for the case of Fig. \ref{fig:Dispersion}

\begin{figure}[hbt]
    \centering
    \includegraphics[width=\linewidth]{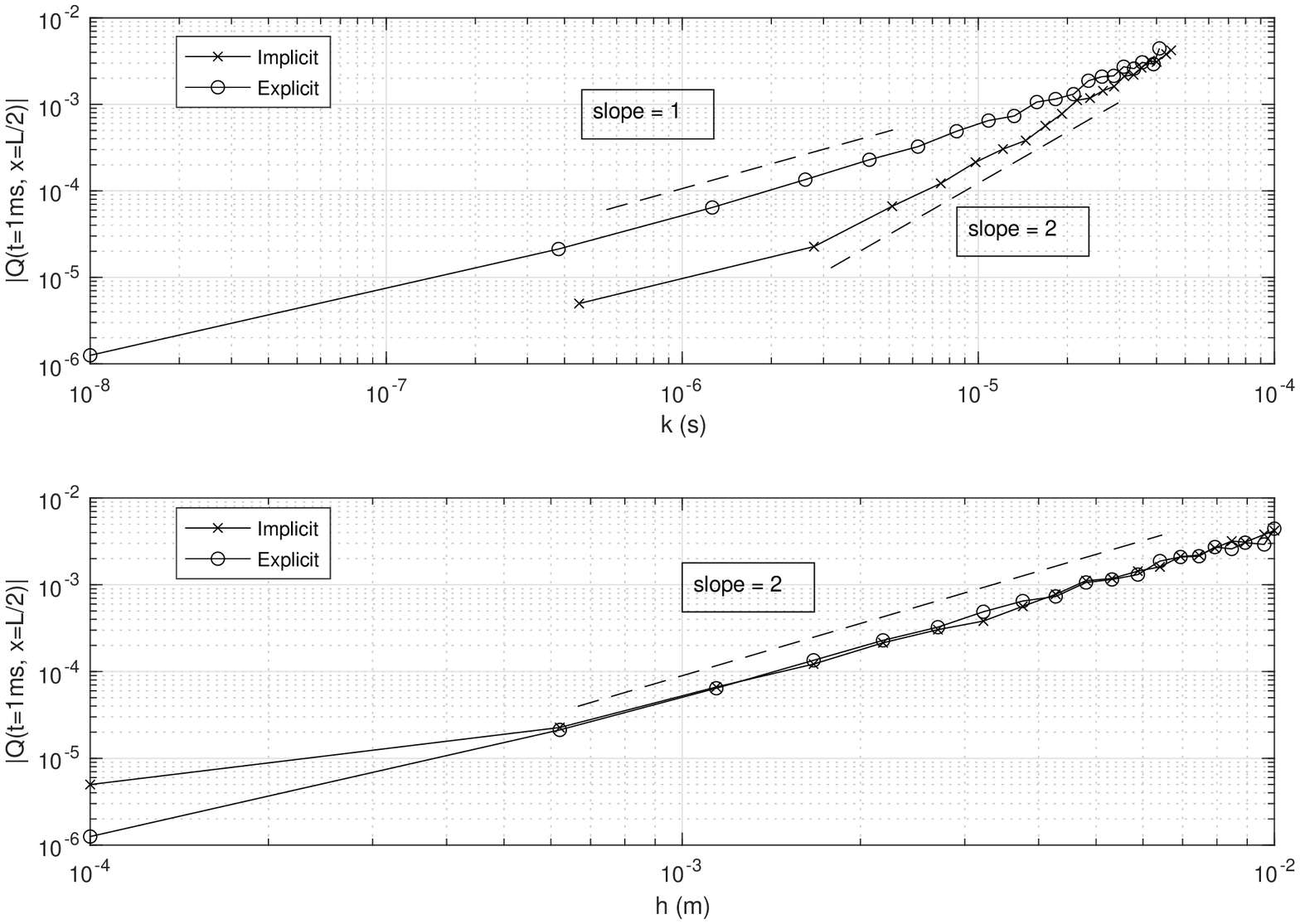}
    \caption{Space-time errors of the implicit and explicit schemes, plotted against $k$ and $h$. Here, the grid spacing is chosen as input parameter for the simulations, and the time step is selected at the corresponding limits of stability for the two schemes. A string with $\rho = 8000$ kg/m$^3$, $E = 2\cdot 10^{11}$ Pa, $r = 0.2$ mm, $L = 1$ m, $T_0=50$ N was used for the simulations. The error is computed as per \eqref{eq:ErrSpaceTime}.}\label{fig:SpaceTimeError}
\end{figure}
While the explicit scheme is formally first-order accurate in time, the discussion in \ref{app:thetaErrs} shows that this fact alone offers a slightly misleading interpretation of the performance of the scheme. This is because the time step and the grid spacing are chosen along the path in the $(h,k)$ plane dictated by the stability condition, and it was shown in \ref{app:thetaErrs} that the scheme is second-order accurate in the path-length $s$. Since $s = O(h)$, the grid spacing can be taken measure of comparison when comparing against the implicit scheme. In practice, it is fair to compare the error curves between the implicit and explicit schemes, {using the same grid spacing, and choosing the time step at the corresponding limits of stability}. This is done if Fig. \ref{fig:SpaceTimeError} below, showing that, as expected, the explicit scheme is first-order in $k$, but the error curves of the two schemes are almost identical when plotted against $h$. Furthermore, time domain simulations in Matlab show once more that the explicit scheme is faster than the explicit one, when the same number of grid points is used. Thus, it can be concluded tha the explicit scheme outperforms the implicit one.
Matlab scripts are available on the companion webpage\footnote{https://mdphys.org/nonlinearString.html}.

\clearpage


\bibliographystyle{elsarticle-num-names}
\bibliography{JSV2020biblio.bib}

\begin{thebibliography}{52}
\providecommand{\natexlab}[1]{#1}
\providecommand{\url}[1]{\texttt{#1}}
\providecommand{\urlprefix}{URL }
\expandafter\ifx\csname urlstyle\endcsname\relax
  \providecommand{\doi}[1]{doi:\discretionary{}{}{}#1}\else
  \providecommand{\doi}[1]{doi:\discretionary{}{}{}\begingroup
  \urlstyle{rm}\url{#1}\endgroup}\fi
\providecommand{\bibinfo}[2]{#2}

\bibitem[{Ruiz(1969)}]{ruiz_thesis}
\bibinfo{author}{P.~Ruiz}, \bibinfo{title}{A technique for simulating the
  vibrations of strings with a digital computer}, Master's thesis,
  \bibinfo{school}{University of Illinois}, \bibinfo{year}{1969}.

\bibitem[{Hiller and Ruiz(1971)}]{ruiz_jaes_1971}
\bibinfo{author}{L.~Hiller}, \bibinfo{author}{P.~Ruiz},
  \bibinfo{title}{Synthesizing Musical Sounds by Solving the Wave Equation for
  Vibrating Objects: Part {I}}, \bibinfo{journal}{{J} {A}udio {E}ng {S}oc}
  \bibinfo{volume}{19}~(\bibinfo{number}{6}) (\bibinfo{year}{1971})
  \bibinfo{pages}{462--470}.

\bibitem[{Bacon and Bowsher(1978)}]{bacon_bowsher_acustica_1978}
\bibinfo{author}{R.~Bacon}, \bibinfo{author}{J.~Bowsher}, \bibinfo{title}{A
  discrete model for of the struck string}, \bibinfo{journal}{Acustica}
  \bibinfo{volume}{41} (\bibinfo{year}{1978}) \bibinfo{pages}{21--27}.

\bibitem[{Chaigne and Askenfelt(1994)}]{chaigne_jasa_1994}
\bibinfo{author}{A.~Chaigne}, \bibinfo{author}{A.~Askenfelt},
  \bibinfo{title}{Numerical simulations of piano strings 1. {A} physical model
  for a struck string using finite difference methods}, \bibinfo{journal}{{J}
  {Acoust} {Soc} {Am}} \bibinfo{volume}{95}~(\bibinfo{number}{2})
  (\bibinfo{year}{1994}) \bibinfo{pages}{1112--1118}.

\bibitem[{Bensa et~al.(2003)Bensa, Bilbao, Kronland-Martinet, and
  Smith}]{bensa_jasa_2003}
\bibinfo{author}{J.~Bensa}, \bibinfo{author}{S.~Bilbao},
  \bibinfo{author}{R.~Kronland-Martinet}, \bibinfo{author}{J.~O. Smith},
  \bibinfo{title}{The simulation of piano string vibration: From physical
  models to finite difference schemes and digital waveguides},
  \bibinfo{journal}{{J} {Acoust} {Soc} {Am}}
  \bibinfo{volume}{114}~(\bibinfo{number}{2}) (\bibinfo{year}{2003})
  \bibinfo{pages}{1095--1107}.

\bibitem[{Ducasse(2005)}]{ducasse_jasa_2005}
\bibinfo{author}{{\'E}.~Ducasse}, \bibinfo{title}{On waveguide modeling of
  stiff piano strings}, \bibinfo{journal}{{J} {Acoust} {Soc} {Am}}
  \bibinfo{volume}{118}~(\bibinfo{number}{3}) (\bibinfo{year}{2005})
  \bibinfo{pages}{1776--1781}.

\bibitem[{Carrier(1945)}]{carrier_QAM_1945}
\bibinfo{author}{G.~F. Carrier}, \bibinfo{title}{On the non-linear vibration
  problem of the elastic string}, \bibinfo{journal}{Quart Appl Math}
  \bibinfo{volume}{3}~(\bibinfo{number}{2}) (\bibinfo{year}{1945})
  \bibinfo{pages}{157--165}.

\bibitem[{Bilbao(2009)}]{bilbao_book}
\bibinfo{author}{S.~Bilbao}, \bibinfo{title}{Numerical {Sound} {Synthesis}:
  {Finite} {Difference} Schemes and Simulation in Musical Acoustics},
  \bibinfo{publisher}{Wiley}, \bibinfo{address}{Chichester, UK},
  \bibinfo{year}{2009}.

\bibitem[{H{\'e}lie and Roze(2008)}]{helie_JSV_2008}
\bibinfo{author}{T.~H{\'e}lie}, \bibinfo{author}{D.~Roze},
  \bibinfo{title}{Sound synthesis of a nonlinear string using Volterra series},
  \bibinfo{journal}{{J} {S}ound {V}ib}
  \bibinfo{volume}{314}~(\bibinfo{number}{1}) (\bibinfo{year}{2008})
  \bibinfo{pages}{275--306}.

\bibitem[{Narasimha(1968)}]{narashima_JSV_1968}
\bibinfo{author}{R.~Narasimha}, \bibinfo{title}{Non-Linear vibration of an
  elastic string}, \bibinfo{journal}{{J} {S}ound {V}ib}
  \bibinfo{volume}{8}~(\bibinfo{number}{1}) (\bibinfo{year}{1968})
  \bibinfo{pages}{134--146}.

\bibitem[{Rowland(2011)}]{rowland_EJP_2011}
\bibinfo{author}{D.~R. Rowland}, \bibinfo{title}{The potential energy density
  in transverse string waves depends critically on longitudinal motion},
  \bibinfo{journal}{Eur {J} {P}hys} \bibinfo{volume}{32}~(\bibinfo{number}{6})
  (\bibinfo{year}{2011}) \bibinfo{pages}{1475--1484}.

\bibitem[{Morse and Ingard(1968)}]{morse_book}
\bibinfo{author}{P.~Morse}, \bibinfo{author}{U.~Ingard},
  \bibinfo{title}{Theoretical Acoustics}, \bibinfo{publisher}{Princeton
  University Press}, \bibinfo{address}{Princeton, NJ, USA},
  \bibinfo{year}{1968}.

\bibitem[{Bilbao(2005)}]{bilbao_jasa_2005}
\bibinfo{author}{S.~Bilbao}, \bibinfo{title}{Conservative numerical methods for
  nonlinear strings}, \bibinfo{journal}{{J} {Acoust} {Soc} {Am}}
  \bibinfo{volume}{118}~(\bibinfo{number}{5}) (\bibinfo{year}{2005})
  \bibinfo{pages}{3316--3327}.

\bibitem[{Ducceschi and Bilbao(2016{\natexlab{a}})}]{ducceschi_st_ICA_2016}
\bibinfo{author}{M.~Ducceschi}, \bibinfo{author}{S.~Bilbao},
  \bibinfo{title}{Modelling Collisions of nonlinear strings against rigid
  barriers: Conservative Finite Difference Schemes With Application to Sound
  Synthesis}, in: \bibinfo{booktitle}{Proc. Int Conf Acous (ICA 2016)},
  \bibinfo{address}{Buenos Aires, Argentina},
  \bibinfo{year}{2016}{\natexlab{a}}.

\bibitem[{Bank and Sujbert(2005)}]{bank_JASA_2005}
\bibinfo{author}{B.~Bank}, \bibinfo{author}{L.~Sujbert},
  \bibinfo{title}{Generation of longitudinal vibrations in piano strings:
  {F}rom physics to sound synthesis}, \bibinfo{journal}{{J} {Acoust} {Soc}
  {Am}} \bibinfo{volume}{117}~(\bibinfo{number}{4}) (\bibinfo{year}{2005})
  \bibinfo{pages}{2268--2278}.

\bibitem[{Ta-Tsien(1994)}]{tatsien_book}
\bibinfo{author}{L.~Ta-Tsien}, \bibinfo{title}{Global classical solutions for
  quasilinear hyperbolic systems}, \bibinfo{publisher}{Wiley},
  \bibinfo{address}{Michigan, USA}, \bibinfo{year}{1994}.

\bibitem[{Chabassier and Joly(2010)}]{chabassier_CMAME_2010}
\bibinfo{author}{J.~Chabassier}, \bibinfo{author}{P.~Joly},
  \bibinfo{title}{Energy preserving schemes for nonlinear {H}amiltonian systems
  of wave equations: {A}pplication to the vibrating piano string},
  \bibinfo{journal}{{C}omput {M}ethod {A}ppl {Mech}}
  \bibinfo{volume}{199}~(\bibinfo{number}{45}) (\bibinfo{year}{2010})
  \bibinfo{pages}{2779--2795}.

\bibitem[{Chabassier(2012)}]{phdthesisChab}
\bibinfo{author}{J.~Chabassier}, \bibinfo{title}{Mod{\'e}lisation et simulation
  num{\'e}rique d'un piano par mod\`{e}les physiques. ({M}odeling and
  simulation of a piano by physical modelling.)}, Ph.D. thesis,
  \bibinfo{school}{Ecole Polytechnique X}, \bibinfo{address}{Paris},
  \bibinfo{year}{2012}.

\bibitem[{Marazzato et~al.(2019)Marazzato, Ern, Mariotti, and
  Monasse}]{marazzato_CMAME_2019}
\bibinfo{author}{F.~Marazzato}, \bibinfo{author}{A.~Ern},
  \bibinfo{author}{C.~Mariotti}, \bibinfo{author}{L.~Monasse},
  \bibinfo{title}{An explicit pseudo-energy conserving time-integration scheme
  for Hamiltonian dynamics}, \bibinfo{journal}{Comput Method Appl Mech}
  \bibinfo{volume}{347} (\bibinfo{year}{2019}) \bibinfo{pages}{906 -- 927}.

\bibitem[{Chatziioannou et~al.(2017)Chatziioannou, Schmutzhard, and
  Bilbao}]{chatziioannou_DAFX_2017}
\bibinfo{author}{V.~Chatziioannou}, \bibinfo{author}{S.~Schmutzhard},
  \bibinfo{author}{S.~Bilbao}, \bibinfo{title}{On Iterative Solutions for
  Numerical Collision Models}, in: \bibinfo{booktitle}{{Proc} {Int} {Conf} {On}
  {Digital} {Audio} {Effects} ({DAFx}-17)}, \bibinfo{address}{Edinburgh, UK},
  \bibinfo{year}{2017}.

\bibitem[{Yang(2016)}]{yang_JCP_2016}
\bibinfo{author}{X.~Yang}, \bibinfo{title}{Linear, first and second-order,
  unconditionally energy stable numerical schemes for the phase field model of
  homopolymer blends}, \bibinfo{journal}{{J} {Comput} {Phys}}
  \bibinfo{volume}{327} (\bibinfo{year}{2016}) \bibinfo{pages}{294--316}, ISSN
  \bibinfo{issn}{0021-9991}.

\bibitem[{Yang(2017)}]{yang_CMAME_2017}
\bibinfo{author}{X.~Yang}, \bibinfo{title}{Linear and unconditionally energy
  stable schemes for the binary fluid-surfactant phase field model},
  \bibinfo{journal}{Comp Methods Appl Mech Eng} \bibinfo{volume}{318}
  (\bibinfo{year}{2017}) \bibinfo{pages}{1005--1029}.

\bibitem[{Zhao et~al.(2017)Zhao, Wang, and Yang}]{zhao_IJNM_2017}
\bibinfo{author}{J.~Zhao}, \bibinfo{author}{Q.~Wang},
  \bibinfo{author}{X.~Yang}, \bibinfo{title}{Numerical approximations for a
  phase field dendritic crystal growth model based on the invariant energy
  quadratization approach}, \bibinfo{journal}{{I}nt {J} {N}umer {M}eth {E}ng}
  \bibinfo{volume}{110}~(\bibinfo{number}{3}) (\bibinfo{year}{2017})
  \bibinfo{pages}{279--300}.

\bibitem[{Chabassier et~al.(2013)Chabassier, Chaigne, and
  Joly}]{chabassier_JASA_2013}
\bibinfo{author}{J.~Chabassier}, \bibinfo{author}{A.~Chaigne},
  \bibinfo{author}{P.~Joly}, \bibinfo{title}{Modeling and simulation of a grand
  piano}, \bibinfo{journal}{{J} {Acoust} {Soc} {Am}}
  \bibinfo{volume}{134}~(\bibinfo{number}{1}) (\bibinfo{year}{2013})
  \bibinfo{pages}{648--665}.

\bibitem[{Han et~al.(1999)Han, Benaroya, and Wei}]{han_JSV_1999}
\bibinfo{author}{S.~Han}, \bibinfo{author}{H.~Benaroya},
  \bibinfo{author}{T.~Wei}, \bibinfo{title}{Dynamics of transversely vibrating
  beams using four engineering theories}, \bibinfo{journal}{{J} {Sound} {Vib}}
  \bibinfo{volume}{225}~(\bibinfo{number}{3}) (\bibinfo{year}{1999})
  \bibinfo{pages}{935--988}.

\bibitem[{Ducceschi and Bilbao(2016{\natexlab{b}})}]{ducceschi_JASA_2016}
\bibinfo{author}{M.~Ducceschi}, \bibinfo{author}{S.~Bilbao},
  \bibinfo{title}{Linear stiff string vibrations in musical acoustics:
  Assessment and comparison of models}, \bibinfo{journal}{{J} {Acoust} {Soc}
  {Am}} \bibinfo{volume}{140}~(\bibinfo{number}{4})
  (\bibinfo{year}{2016}{\natexlab{b}}) \bibinfo{pages}{2445--2454}.

\bibitem[{Lopes et~al.(2015)Lopes, H{\'e}lie, and Falaize}]{lopes_IFAC_2015}
\bibinfo{author}{N.~Lopes}, \bibinfo{author}{T.~H{\'e}lie},
  \bibinfo{author}{A.~Falaize}, \bibinfo{title}{Explicit second-order accurate
  method for the passive guaranteed simulation of port-Hamiltonian systems},
  in: \bibinfo{booktitle}{{Proc} {5th IFAC} 2015}, \bibinfo{address}{Lyon,
  France}, \bibinfo{year}{2015}.

\bibitem[{Falaize(2016)}]{phdthesis}
\bibinfo{author}{A.~Falaize}, \bibinfo{title}{Mod{\'e}lisation, simulation,
  g{\'e}n{\'e}ration de code et correction de syst{\`e}mes multi-physiques
  audios: Approche par r{\'e}seau de composants et formulation {H}amiltonienne
  {\`A} Ports. ({M}odeling, Simulation, code generation and correction of
  multiphysics audio systems: Component network approach and
  {P}ort-{H}amiltonian formulation)}, Ph.D. thesis,
  \bibinfo{school}{Universit{\'e} Pierre et Marie Curie},
  \bibinfo{address}{Paris}, \bibinfo{year}{2016}.

\bibitem[{Falaize and H{\'e}lie(2016)}]{falaize:hal-01390501}
\bibinfo{author}{A.~Falaize}, \bibinfo{author}{T.~H{\'e}lie},
  \bibinfo{title}{{Passive Guaranteed Simulation of Analog Audio Circuits: A
  Port-Hamiltonian Approach}}, \bibinfo{journal}{{Appl} Sci}
  \bibinfo{volume}{6} (\bibinfo{year}{2016}) \bibinfo{pages}{273--273}.

\bibitem[{Lopes(2016)}]{Lopesphdthesis}
\bibinfo{author}{N.~Lopes}, \bibinfo{title}{Approche passive pour la
  mod{\'e}lisation, la simulation et l’{\'e}tude d’un banc de test
  robotis{\'e} pour les instruments de type cuivre. (Passive approach for
  modelling, simulation and study of a robotic test bench for brass
  instruments.)}, Ph.D. thesis, \bibinfo{school}{Universit{\'e} Pierre et Marie
  Curie}, \bibinfo{address}{Paris}, \bibinfo{year}{2016}.

\bibitem[{Jiang. et~al.(2019)Jiang., Cai, and Wang}]{jiang_JSC_2019}
\bibinfo{author}{C.~Jiang.}, \bibinfo{author}{W.~Cai},
  \bibinfo{author}{Y.~Wang}, \bibinfo{title}{A Linearly Implicit and Local
  Energy-Preserving Scheme for the Sine-Gordon Equation Based on the Invariant
  Energy Quadratization Approach}, \bibinfo{journal}{{J} {S}ci {C}omput}
  \bibinfo{volume}{80} (\bibinfo{year}{2019}) \bibinfo{pages}{1629--1655}.

\bibitem[{Ducceschi and Bilbao(2019{\natexlab{a}})}]{ducceschi_ICA_2019}
\bibinfo{author}{M.~Ducceschi}, \bibinfo{author}{S.~Bilbao},
  \bibinfo{title}{Non-iterative, conservative schemes for geometrically exact
  nonlinear string vibration}, in: \bibinfo{booktitle}{Proc {I}nt {C}onf
  {A}coust (ICA 2019)}, \bibinfo{address}{Aachen, Germany},
  \bibinfo{year}{2019}{\natexlab{a}}.

\bibitem[{Ducceschi et~al.(2021{\natexlab{a}})Ducceschi, Bilbao, Willemsen, and
  Serafin}]{ducceschi_jasa_2021}
\bibinfo{author}{M.~Ducceschi}, \bibinfo{author}{S.~Bilbao},
  \bibinfo{author}{S.~Willemsen}, \bibinfo{author}{S.~Serafin},
  \bibinfo{title}{Linearly-implicit schemes for collisions in musical acoustics
  based on energy quadratisation}, \bibinfo{journal}{J. Acoust. Soc. Am.}
  \bibinfo{volume}{149}~(\bibinfo{number}{5})
  (\bibinfo{year}{2021}{\natexlab{a}}) \bibinfo{pages}{3502--3516}.

\bibitem[{Torin and Bilbao(2013)}]{torin_DAFX_2013}
\bibinfo{author}{A.~Torin}, \bibinfo{author}{S.~Bilbao}, \bibinfo{title}{A 3D
  multi-plate environment for sound synthesis}, in: \bibinfo{booktitle}{{Proc}
  {Int} {Conf} {On} {Dig} {Audio} {Eff} ({DAFx}-13)},
  \bibinfo{address}{Maynooth, Ireland}, \bibinfo{year}{2013}.

\bibitem[{Gazdag(1973)}]{gazdag_JCP_1973}
\bibinfo{author}{J.~Gazdag}, \bibinfo{title}{Numerical convective schemes based
  on accurate computation of space derivatives}, \bibinfo{journal}{{J} {Comput}
  {Phys}} \bibinfo{volume}{13}~(\bibinfo{number}{1}) (\bibinfo{year}{1973})
  \bibinfo{pages}{100 -- 113}.

\bibitem[{Gazdag(1981)}]{gazdag_geophysics_1981}
\bibinfo{author}{J.~Gazdag}, \bibinfo{title}{Modeling of the acoustic wave
  equation with transform methods}, \bibinfo{journal}{Geophysics}
  \bibinfo{volume}{46}~(\bibinfo{number}{6}) (\bibinfo{year}{1981})
  \bibinfo{pages}{854--859}.

\bibitem[{Ducceschi and Bilbao(2019{\natexlab{b}})}]{ducceschi_WaveMotion_2019}
\bibinfo{author}{M.~Ducceschi}, \bibinfo{author}{S.~Bilbao},
  \bibinfo{title}{Conservative finite difference time domain schemes for the
  prestressed {T}imoshenko, shear and {E}uler-{B}ernoulli beam equations},
  \bibinfo{journal}{Wave Motion} \bibinfo{volume}{89}
  (\bibinfo{year}{2019}{\natexlab{b}}) \bibinfo{pages}{142 -- 165}.

\bibitem[{LeVeque(2007)}]{leveque_book}
\bibinfo{author}{R.~LeVeque}, \bibinfo{title}{Finite Difference Methods for
  Ordinary and Partial Differential Equations. Steady State and Time Dependent
  Problems}, \bibinfo{publisher}{SIAM}, \bibinfo{address}{Philadelphia, USA},
  \bibinfo{year}{2007}.

\bibitem[{Joly(2003)}]{joly_book}
\bibinfo{author}{P.~Joly}, \bibinfo{title}{Variational methods for
  time-dependent wave propagation problems}, in:
  \bibinfo{editor}{J.~Fagerberg}, \bibinfo{editor}{D.~C. Mowery},
  \bibinfo{editor}{R.~R. Nelson} (Eds.), \bibinfo{booktitle}{Topics in
  Computational Wave Propagation, Lecture Notes in Computational Science and
  Engineering}, \bibinfo{publisher}{Springer}, \bibinfo{address}{Berlin},
  \bibinfo{pages}{201--264}, \bibinfo{year}{2003}.

\bibitem[{Chabassier and Imperiale(2017)}]{chabassier_CRM_2017}
\bibinfo{author}{J.~Chabassier}, \bibinfo{author}{S.~Imperiale},
  \bibinfo{title}{Space/time convergence analysis of a class of conservative
  schemes for linear wave equations}, \bibinfo{journal}{Comptes Rendus
  Mathematique} \bibinfo{volume}{355}~(\bibinfo{number}{3})
  (\bibinfo{year}{2017}) \bibinfo{pages}{282--289}.

\bibitem[{Cohen and Joly(1996)}]{cohen_1996}
\bibinfo{author}{G.~Cohen}, \bibinfo{author}{P.~Joly},
  \bibinfo{title}{Construction analysis of fourth-order finite difference
  schemes for the acoustic wave equation in nonhomogeneous media},
  \bibinfo{journal}{SIAM {J} {Numer} {Anal}}
  \bibinfo{volume}{33}~(\bibinfo{number}{4}) (\bibinfo{year}{1996})
  \bibinfo{pages}{1266--1302}.

\bibitem[{Joly and Rodr{\'i}guez(2010)}]{joly_2010}
\bibinfo{author}{P.~Joly}, \bibinfo{author}{J.~Rodr{\'i}guez},
  \bibinfo{title}{Optimized higher order time discretization of second order
  hyperbolic problems: Construction and numerical study}, \bibinfo{journal}{{J}
  {Comput} {Appl} {Math}} \bibinfo{volume}{234}~(\bibinfo{number}{6})
  (\bibinfo{year}{2010}) \bibinfo{pages}{1953--1961}.

\bibitem[{Dablian(1986)}]{dablian_2010}
\bibinfo{author}{M.~Dablian}, \bibinfo{title}{The application of high order
  differencing for the scalar wave equation}, \bibinfo{journal}{Geophysics}
  \bibinfo{volume}{51} (\bibinfo{year}{1986}) \bibinfo{pages}{54--56}.

\bibitem[{Warming and Hyett(1974)}]{warming_1974}
\bibinfo{author}{R.~Warming}, \bibinfo{author}{B.~Hyett}, \bibinfo{title}{The
  modified equation approach to the stability and accuracy analysis of
  finite-difference methods}, \bibinfo{journal}{{J} {Comput} {Phys}}
  \bibinfo{volume}{14}~(\bibinfo{number}{2}) (\bibinfo{year}{1974})
  \bibinfo{pages}{159--179}.

\bibitem[{Bilbao and Hamilton(2018)}]{bilbao_2018}
\bibinfo{author}{S.~Bilbao}, \bibinfo{author}{B.~Hamilton},
  \bibinfo{title}{Higher-order accurate two-step finite difference schemes for
  the many-dimensional wave equation}, \bibinfo{journal}{{J} {Comput} {Phys}}
  \bibinfo{volume}{367} (\bibinfo{year}{2018}) \bibinfo{pages}{134--165}.

\bibitem[{Germain and Werner(2015)}]{germain_dafx_2015}
\bibinfo{author}{F.~G. Germain}, \bibinfo{author}{K.~J. Werner},
  \bibinfo{title}{Design principles for lumped model discretisation using
  Moebius transforms}, in: \bibinfo{booktitle}{Proc Digital Audio Effects
  (DAFx-15)}, \bibinfo{address}{Trondheim, Norway}, \bibinfo{year}{2015}.

\bibitem[{Ducceschi et~al.(2021{\natexlab{b}})Ducceschi, Bilbao, and
  Webb}]{ducceschi_DAFX_2021}
\bibinfo{author}{M.~Ducceschi}, \bibinfo{author}{S.~Bilbao},
  \bibinfo{author}{C.~Webb}, \bibinfo{title}{Non-iterative Schemes For The
  Simulation Of Nonlinear Audio Circuits}, in: \bibinfo{booktitle}{{Proc}
  {Digital} {Audio} {Effects} ({DAFx}-21)}, \bibinfo{address}{Vienna, Austria},
  \bibinfo{year}{2021}{\natexlab{b}}.

\bibitem[{Lele(1992)}]{lele_1992}
\bibinfo{author}{S.~K. Lele}, \bibinfo{title}{Compact finite difference schemes
  with spectral-like resolution}, \bibinfo{journal}{J {Comput} {Phys}}
  \bibinfo{volume}{103}~(\bibinfo{number}{1}) (\bibinfo{year}{1992})
  \bibinfo{pages}{16--42}.

\bibitem[{Bilbao et~al.(2015)Bilbao, Torin, and
  Chatziioannou}]{bilbao_AAUA_2015}
\bibinfo{author}{S.~Bilbao}, \bibinfo{author}{A.~Torin},
  \bibinfo{author}{V.~Chatziioannou}, \bibinfo{title}{Numerical Modeling of
  Collisions in Musical Instruments}, \bibinfo{journal}{Acta {A}cust {U}nited
  with {A}cust} \bibinfo{volume}{101} (\bibinfo{year}{2015})
  \bibinfo{pages}{155--173}.

\bibitem[{Ducceschi and Bilbao(2019{\natexlab{c}})}]{ducceschi_DAFx_2019}
\bibinfo{author}{M.~Ducceschi}, \bibinfo{author}{S.~Bilbao},
  \bibinfo{title}{Non-iterative solvers for nonlinear problems: the case of
  collisions}, in: \bibinfo{booktitle}{Proc 22nd Conf of Digital Audio Effects
  (DAFx-19)}, \bibinfo{address}{Birmingham, UK},
  \bibinfo{year}{2019}{\natexlab{c}}.

\bibitem[{Chatziioannou and van Walstijn(2015)}]{chatziioannou_JSV_2015}
\bibinfo{author}{V.~Chatziioannou}, \bibinfo{author}{M.~van Walstijn},
  \bibinfo{title}{Energy conserving schemes for the simulation of musical
  instrument contact dynamics}, \bibinfo{journal}{{J} {Sound} {Vib}}
  \bibinfo{volume}{339} (\bibinfo{year}{2015}) \bibinfo{pages}{262--279}.

\bibitem[{Ducceschi et~al.(2014)Ducceschi, Cadot, Touz\'e, and
  Bilbao}]{ducceschi_PhysD_2014}
\bibinfo{author}{M.~Ducceschi}, \bibinfo{author}{O.~Cadot},
  \bibinfo{author}{C.~Touz\'e}, \bibinfo{author}{S.~Bilbao},
  \bibinfo{title}{Dynamics of the wave turbulence spectrum in vibrating plates:
  A numerical investigation using a conservative finite difference scheme},
  \bibinfo{journal}{Physica D} \bibinfo{volume}{280-281} (\bibinfo{year}{2014})
  \bibinfo{pages}{73--85}.

\end{thebibliography}







\end{document}